\documentclass[11pt]{article}
\usepackage{amsfonts}
\usepackage{color}
\usepackage{amsmath,amssymb,comment}
\newtheorem{theo}{Theorem}[section]
\newtheorem{lem}[theo]{Lemma}
\newtheorem{cor}[theo]{Corollary}

\newtheorem{defi}[theo]{Definition}

\newcommand{\mysection}[1]{\section{#1} \setcounter{equation}{0}}
\newcommand{\proof}{{\sc Proof.} \quad}
\newcommand{\proofc}{{\sc Proof} \ }
\newcommand{\be}{\begin{equation} \label}
\newcommand{\ee}{\end{equation}}
\newcommand{\bea}{\begin{eqnarray}\label}
\newcommand{\eea}{\end{eqnarray}}
\newcommand{\bas}{\begin{eqnarray*}}
\newcommand{\eas}{\end{eqnarray*}}
\newcommand{\bit}{\begin{itemize}}
\newcommand{\eit}{\end{itemize}}
\newcommand{\qed}{\hfill$\Box$ \vskip.2cm}
\newcommand{\nn}{\nonumber}
\newcommand{\R}{\mathbb{R}}
\newcommand{\N}{\mathbb{N}}
\newcommand{\pO}{\partial\Omega}

\newcommand{\eps}{\varepsilon}

\newcommand{\wto}{\rightharpoonup}
\newcommand{\wsto}{\stackrel{\star}{\rightharpoonup}}
\newcommand{\hra}{\hookrightarrow}
\newcommand{\io}{\int_\Omega}
\newcommand{\na}{\nabla}
\newcommand{\Del}{\Delta}

\newcommand{\al}{\alpha}

\newcommand{\pa}{\partial}
\newcommand{\bom}{\overline{\Omega}}
\newcommand{\Om}{\Omega}

\newcommand{\ov}{\overline}

\newcommand{\wt}{\widetilde}

\newcommand{\hs}{\hspace*}
\newcommand{\sm}{\setminus}

\newcommand{\vp}{\varphi}
\newcommand{\lbal}{\left\{ \begin{array}{l}}
\newcommand{\lball}{\left\{ \begin{array}{ll}}
\newcommand{\ear}{\end{array} \right.}
\newcommand{\ouz}{\ov{u}_0}

\newcommand{\abs}{\\[5pt]}

\newcommand{\adb}{\allowdisplaybreaks}

\newcommand{\ueps}{u_\eps}
\newcommand{\veps}{v_\eps}

\newcommand{\pe}{\phi_\eps}
\newcommand{\epss}{\eps_\star}
\newcommand{\xis}{\xi_\star}
\newcommand{\Ieps}{I_\eps}
\begin{document}
\adb
%
%
\title{A strongly degenerate migration-consumption model\\
in domains of arbitrary dimension}
\author{
Michael Winkler\footnote{michael.winkler@math.uni-paderborn.de}\\
{\small Universit\"at Paderborn, Institut f\"ur Mathematik}\\
{\small 33098 Paderborn, Germany} }
\date{}
\maketitle
\begin{abstract}
\noindent 
In a smoothly bounded convex domain $\Omega\subset\R^n$ with $n\ge 1$,
a no-flux initial-boundary value problem for
\bas
	\lbal
	u_t=\Del \big(u\phi(v)\big), \\[1mm]
	v_t=\Del v-uv,
	\ear
\eas
is considered
under the assumption that near the origin, the function $\phi$ suitably generalizes the prototype given by
\bas
	\phi(\xi)=\xi^\al,
	\qquad \xi\in [0,\xi_0].
\eas
By means of separate approaches, it is shown that in both cases $\al\in (0,1)$ and $\al\in [1,2]$ some global
weak solutions exist which, inter alia, satisfy
\bas
	C(T):= {\rm ess} \sup_{\hs{-5mm} t\in (0,T)} \io u(\cdot,t)\ln u(\cdot,t) < \infty
	\qquad \mbox{for all } T>0,
\eas
with $\sup_{T>0} C(T)<\infty$ if $\al\in [1,2]$.\abs
\end{abstract}
\newpage
\section{Introduction}\label{intro}
Local sensing mechanisms are relevant to partially directed motion of cell motion (\cite{maini_et_al_JTB2013}, \cite{liu}, 
\cite{fu}, \cite{othmer_stevens}, \cite{engwer_hunt_surulescu}).
It their simplest form, macroscopic models for such processes 
describe the population density $u=u(x,t)$ by parabolic equations of the form
\be{00}
	u_t=\Del \big( a(x,t) u \big)
\ee
(\cite{maini_et_al_JTB2013}, \cite{engwer_hunt_surulescu}), where in typical application situations,
the cell motility coefficient $a$ may depend on a chemical substance, represented through its concentration $v=v(x,t)$
itself forming an unknown of the system, via various functional laws (\cite{fu}, \cite{KS3}).\abs
In recent analytical literature, significant activity has been directed toward an understanding of resulting two-component
parabolic models in cases in which the respective signal is produced by cells, and which thus, as in the classical 
Keller-Segel systems from \cite{KS1}, reflect taxis-mediated active communication between cells.
Hence focusing on systems of the form
\be{01}
	\lbal
	u_t=\Del \big(u\phi(v)\big), \\[1mm]
	v_t=\Del v-v+u,
	\ear
\ee
as well as on some parabolic-elliptic simplifications thereof, a considerable collection of studies has identified 
various conditions on the key ingredient $\phi$ as sufficient for global solvability and hence for suppression
of finite-time blow-up (see \cite{ahn_yoon}, \cite{fujie_jiang_arxiv}, \cite{fujie_jiang_JDE2020}, 
\cite{fujie_jiang_ACAP2021}, \cite{jiang_laurencot},
\cite{burger}, \cite{fujie_senba}, \cite{taowin_M3AS} for an incomplete selection, and also
\cite{jin_kim_wang}, \cite{wenbin_lv_ZAMP}, \cite{wenbin_lv_EECT}, \cite{wenbin_lv_PROCA}, \cite{wang_wang}, \cite{yifu_wang} and 
\cite{win_NON} for some studies on variants accounting for sources and density-dependent diffusion mechanisms);
on the other hand, some results detecting the occurrence of infinite-time blow-up in the particular case when 
$\phi(v)=e^{-v}$, $v\ge 0$, indicate a certain reminiscence of Keller-Segel dynamics (\cite{fujie_jiang_arxiv}; cf.also
\cite{fujie_senba}):
Passing over from classical Keller-Segel-production systems to models of the form (\ref{01}) may thus, depending on the 
choice of $\phi$, delay but not entirely rule out unboundedness phenomena (\cite{burger}).\abs
In comparison to the above, much less seems known for related systems addressing situations of local sensing
in which the directing signal is
consumed by individuals, and in which thus cells, in particular, are incapable of active communication.
In fact, for corresponding migration-absorption models of the form
\be{02}
	\lbal
	u_t=\Del \big(u\phi(v)\big), \\[1mm]
	v_t=\Del v-uv,
	\ear
\ee
the literature so far appears to concentrate on non-degenerate cases determined by motilities which are strictly 
positive on $[0,\infty)$, and which hence reflect non-degenerate diffusion:
In such situations, the additional dissipative influence exerted by the absorptive reaction substantially facilitates
global existence theories, in frameworks both of classical small-data and of generalized large-data solutions 
(\cite{li_zhao_ZAMP}, \cite{liwin2}); as strongly indicated by quite far-reaching findings on large time stabilization
toward spatially homogeneous steady states, however, non-degenerate settings of this flavor seem unable to adequately
capture any of the strongly structure-supporting features of collective movement observed in populations
of aerobic bacteria (\cite{kawasaki_JTB1997}).\abs
This is in line with refined modeling approaches which, in order to particularly address such situations,
suggest to explicitly
account for reduction of bacterial motility in nutrient-poor environments (\cite{kawasaki_JTB1997}, \cite{plaza}).
Indeed, a recent result indicates that in sharp contrast to said case of positive $\phi$,
nontrivial long term dynamics may indeed occur in 
corresponding versions of (\ref{02}) which accordingly include migration rates reflecting 
motility degeneracies at small signal concentrations: 
When $\phi$ is suitably smooth with 
\be{03}
	\phi(0)=0, \quad \phi'(0)>0
	\quad \mbox{and} \quad
	\phi>0 \mbox{ on } (0,\infty),
\ee
namely,
an associated no-flux type boundary value problem for (\ref{02}) in one- or two-dimensional domains has been found to admit
some classical solutions $(u,v)$ for which $u$ approaches a nonconstant profile in the large time limit
(\cite{win_sig_dep_mot_SMP}).
However, the question whether such types of behavior are restricted to such special settings, or rather constitute
a characteristic feature of degenerate motilities in (\ref{02}) within a more general framework, appears to be open up to now;
in particular, the only precedent we are aware of which addresses somewhat stronger degeneracies,
in fact covering any decay behavior of $\phi$ near the origin which is of essentially algebraic type,
is still limited to domains in $\R^n$ with $n\le 2$ (\cite{win_2d}).\abs
{\bf Main results.} \quad
The present manuscript attempts to design an analytical approach for (\ref{02}) which does not only
allow for the inclusion of motility degeneracies more general than those determined by (\ref{03}), but which moreover
does not rely on assumptions on low dimensionality. 
Due to challenges which in comparison to those encountered in the setup from (\ref{03}) seem considerably increased,
we will focus here on issues from basic solution and regularity theories, leaving more detailed qualitative
investigation for future research.\abs
Specifically, we shall consider
the initial-boundary value problem
\be{0}
	\left\{ \begin{array}{ll}	
	u_t = \Del \big(u\phi(v)\big),
	\qquad & x\in\Om, \ t>0, \\[1mm]
	v_t = \Delta v-uv,
	\qquad & x\in\Om, \ t>0, \\[1mm]
	\na \big(u\phi(v)\big) \cdot \nu = \na v \cdot \nu=0,
	\qquad & x\in\pO, \ t>0, \\[1mm]
	u(x,0)=u_0(x), \quad v(x,0)=v_0(x),
	\qquad & x\in\Om,
	\end{array} \right.
\ee
in $n$-dimensional smoothly bounded convex domains $\Om$, under the assumption that near the origin, 
$\phi$ suitably generalizes the prototype given by
\be{proto}
	\phi(\xi)=\xi^\al,
	\qquad \xi\in [0,\xi_0],
\ee
with certain $\al>0$ and $\xi_0>0$.\abs
Our first step will concentrate on the case $\al\in (0,1)$, in which on the one hand a comparatively mild degeneracy
retains some strength of diffusive smoothing, but for which on the other the corresponding cross-diffusive action 
is considerably singular in regions where $v$ is small. In the context of the identity
\be{log}
	\frac{d}{dt} \io u\ln u + \io \phi(v) \frac{|\na u|^2}{u} = - \io \phi'(v) \na u\cdot\na v
\ee
formally determining the evolution of the associated logarithmic entropy, the latter becomes manifest in a singular 
factor $\phi'(v)$ appearing in the rightmost integral, and a straightforward estimation thereof in terms of the 
dissipated quantity on the left thus seems not expedient. 
Forming a key observation in this regard, it will turn out that by linearly combining (\ref{log}) with a corresponding
identity describing the evolution of an appropriate bounded quantity, this expression can be suitably diminished
in strength, and hence become conveniently controllable by respective diffusive contributions.
Indeed, in Section \ref{sect3} we shall see that for appropriately chosen $a>0$, an inequality of the form
\be{energy1}
	\frac{d}{dt} \bigg\{ \io u\ln u - a \io u\phi(v) \bigg\}
	+ \frac{1}{2} \io \phi(v) \frac{|\na u|^2}{u}
	\le C \io \frac{|\na v|^4}{v^3} + C \io u^2 v^\al
\ee
holds for solutions to certain regularized variants of (\ref{0}) (see (\ref{0eps}) and Lemma \ref{lem12}).
In conjunction with some basic regularity features, originating from a standard duality-based reasoning
and providing bounds for both summands on the right of (\ref{energy1}) (Lemma \ref{lem4}, Lemma \ref{lem56} 
and Corollary \ref{cor58}),
this will lead to a priori information sufficient for the derivation of a result on global existence of weak
solutions with locally bounded logarithmic entropies.\abs
More precisely, the first of our main results can be stated as follows.
\begin{theo}\label{theo15}
  Let $n\ge 1$ and $\Omega\subset\R^n$ be a bounded convex domain with smooth boundary, and suppose that 
  \be{preg}
	\phi\in C^0([0,\infty)) \cap C^3((0,\infty))
	\mbox{\quad is such that \quad $\phi(0)=0$ \quad and \quad $\phi(\xi)>0$ \ for all } \xi>0,
  \ee
  and that with some $\al\in (0,1)$ and $\xi_0>0$ we have
  \be{pa}
	\liminf_{\xi\searrow 0} \frac{\phi(\xi)}{\xi^\al}>0
  \ee
  and
  \be{p2}
	\limsup_{\xi\searrow 0} \frac{|\phi'(\xi)|}{\xi^{\al-1}}<\infty
  \ee
  as well as
  \be{p3}
	\big(\phi^\frac{1}{\al}\big)''(\xi) \le 0
	\qquad \mbox{for all } \xi\in (0,\xi_0).
  \ee
  Then whenever
  \be{init}
	\lbal
	u_0\in W^{1,\infty}(\Om) \mbox{ is nonnegative with $u_0\not\equiv 0$}
	\qquad \mbox{and } \\[1mm]
	v_0\in W^{1,\infty}(\Omega) \mbox{ satisfies $v_0>0$ in $\bom$,}
	\ear
  \ee
  one can find
  \be{reg}
	\lbal
	u\in L^\infty((0,\infty);L^1(\Om)) 
	\qquad \mbox{and} \\[1mm]
	v\in L^\infty(\Om\times (0,\infty)) \cap L^\infty_{loc}([0,\infty);W^{1,2}(\Om)) 
		\cap L^2_{loc}([0,\infty);W^{2,2}(\Om))
	\ear
  \ee
  such that $u\ge 0$ and $v>0$ a.e.~in $\Om\times (0,\infty)$, that 
  \be{massu}
	\io u(\cdot,t)=\io u_0
	\qquad \mbox{for a.e.~ } t>0
  \ee
  as well as
  \be{15.2}
	{\rm ess} \sup_{\hs{-5mm} t\in (0,T)} \io u(\cdot,t)\ln u(\cdot,t) <\infty
	\qquad \mbox{for all } T>0
  \ee
  and
  \be{15.3}
	\int_0^T \io \big| \na (u\phi(v))\big|^{p(\al)} <\infty
	\qquad \mbox{for all } T>0
  \ee
  with
  \be{15.4}
	p(\al):=\lball
	\frac{2}{2-\al}
	\qquad & \mbox{if } \al\in (0,\frac{1}{2}), \\[2mm]
	\frac{4}{3}
	\qquad & \mbox{if } \al\in [\frac{1}{2},1),
	\ear
  \ee
  and that $(u,v)$ forms a global weak solution of (\ref{0}) in the sense of Definition \ref{dw}.
\end{theo}
Next concerned with stronger degeneracies, we will need to appropriately cope with an apparent lack of 
structural features comparable to that from (\ref{energy1}) in the case when $\al\ge 1$.
We shall alternatively build our analysis in this respect on the observation that as long as $\al\in (1,2)$,
a favorable entropy-like evolution property of Dirichlet integrals involving mildly singular weights
will lead to an inequality of the form
\bas
	\frac{d}{dt} \io v^{\al-2} |\na v|^2 + \frac{1}{C} \io v^{\al-4} |\na v|^4
	\le C \io u^2 v^\al
\eas
(Lemma \ref{lem6}). Along with fairly standard extensions thereof to the borderline cases $\al=1$ and $\al=2$
(Corollary \ref{cor57}), this can be combined with an again duality-based control
of $\int_0^T \io u^2 v^\al$ which now can even be achieved with bounds uniform with respect to $T$ thanks to the fact that
$uv^\al$ is esentially dominated by the quantity $uv$ known to belong to $L^1(\Om\times (0,\infty))$ according to
the second equation in (\ref{0}) (Lemma \ref{lem44} and (\ref{uv})). 
In Lemma \ref{lem7}, we thereby see that whenever $\al\in [1,2]$, a fairly straightforward estimation of the right-hand side
in (\ref{log}) becomes possible so as to ensure bounds which are now even independent of time.\abs
In conclusion, this will enable us to derive the second of our main results, asserting 
global solvability and temporally uniform bounds in the presence of such superlinear degeneracies, and in 
arbitrarily high-dimensional settings, in the following sense.
\begin{theo}\label{theo155}
  Suppose that $n\ge 1$ and $\Omega\subset\R^n$ is a bounded convex domain with smooth boundary, and that
  $\phi$ satisfies (\ref{preg}), (\ref{pa}) and (\ref{p2}) with some $\al\in [1,2]$.
  Then for any choice of $(u_0,v_0)$ fulfilling (\ref{init}),
  one can find functions $u$ and $v$ which satisfy (\ref{reg}) with
  $u\ge 0$ and $v>0$ a.e.~in $\Om\times (0,\infty)$, which are such that (\ref{massu}) holds as well as
  \be{155.1}
	{\rm ess} \sup_{\hs{-5mm} t>0} \io u(\cdot,t)\ln u(\cdot,t)<\infty,
	\qquad
	{\rm ess} \sup_{\hs{-5mm} t>0} \io |\na v(\cdot,t)|^2 <\infty
  \ee
  and
  \be{155.2}
	\int_0^\infty \io v^{\al-4} |\na v|^4 + \int_0^\infty \io \big| \na (u\phi(v))\big|^\frac{4}{3} <\infty,
  \ee
  and that $(u,v)$ forms a global weak solution of (\ref{0}) in the sense of Definition \ref{dw}.
\end{theo}
We remark that the theory developed here forms the basis of the refined qualitative analysis undertaken in \cite{win_NON2023};
further information on large time stabilization has been obtained in \cite{laurencot_CMS} and in \cite{win_sig_dep_mot_SMP}.
\mysection{Approximation and some basic regularity features}
The solution concept to be considered below appears to be quite in line with standard notions of generalized
solvability in second order parabolic problems, particularly involving one-step integration by parts only.
\begin{defi}\label{dw}
  Suppose that $\phi\in C^0([0,\infty))$, $u_0\in L^1(\Om)$ and $v_0\in L^\infty(\Om)$ are all nonnegative.
  Then by a {\em global weak solution} of (\ref{0}) we mean a pair $(u,v)$ of nonnegative functions
  \be{w1}
	\left\{ \begin{array}{l}
	u\in L^1_{loc}(\bom\times [0,\infty))
	\qquad \mbox{and} \\[1mm]
	v\in L^\infty_{loc}(\bom\times [0,\infty)) \cap L^1_{loc}([0,\infty);W^{1,1}(\Om))
	\ear
  \ee
  such that
  \be{w2}
	\na (u\phi(v)) \in L^1_{loc}(\bom\times [0,\infty);\R^n),
  \ee
  and that
  \bea{wu}
	\int_0^\infty \io u\vp_t
	+ \io u_0 \vp(\cdot,0)
	= \int_0^\infty \io \na (u\phi(v))\cdot\na\vp
  \eea
  and
  \be{wv}
	\int_0^\infty \io v\vp_t  
	+ \io v_0\vp(\cdot,0)
	= \int_0^\infty \io \na v\cdot\na\vp
	+ \int_0^\infty \io uv\vp
  \ee
  hold for any $\vp\in C_0^\infty(\bom\times [0,\infty))$.
\end{defi}
In order to achieve a convenient approximation of (\ref{0}), let us regularize not only the diffusive contribution to the first
equation, but also the reaction part in the second. In fact, this will ensure that for each $\eps\in (0,1)$, the problem
\be{0eps}
	\lball
	u_{\eps t} = \Del \big(\ueps\phi_\eps(\veps)\big),
	\qquad & x\in\Omega, \ t>0, \\[1mm]
	v_{\eps t} = \Delta\veps - \frac{\ueps\veps}{1+\eps\ueps},
	\qquad & x\in\Omega, \ t>0, \\[1mm]
	\frac{\partial\ueps}{\partial\nu}=\frac{\partial\veps}{\partial\nu}=0,
	\qquad & x\in\pO, \ t>0, \\[1mm]
	\ueps(x,0)=u_0(x), \quad
	\veps(x,0)=v_0(x),
	\qquad & x\in\Omega.
	\ear
\ee
with
\be{pe}
	\pe(\xi):=\phi(\xi)+\eps,
	\qquad \xi\ge 0, \ \eps\in (0,1),
\ee
is globally solvable in the classical sense:
\begin{lem}\label{lem_loc}
  Let $n\ge 1$ and $\Om\subset\R^n$ be a bounded domain with smooth boundary, and assume (\ref{preg}) and (\ref{init}).
  Then for each $\eps\in (0,1)$ there exist 
  \be{01.1}
	\lbal
	\ueps\in C^0(\bom\times [0,\infty)) \cap C^{2,1}(\bom\times (0,\infty))
	\qquad \mbox{and} \\[1mm]
	\veps \in \bigcap_{q>n} C^0([0,\infty);W^{1,q}(\Om)) \cap C^{2,1}(\bom\times (0,\infty))
	\ear
  \ee
  such that $\ueps \ge 0$ and $\veps>0$ in $\bom\times [0,\infty)$, 
  that $(\ueps,\veps)$ solves (\ref{0eps})-(\ref{pe}) in the classical sense, and that
  \be{mass}
	\io \ueps(\cdot,t) = \io u_0
	\qquad \mbox{for all } t>0
  \ee
  and
  \be{vinfty}
	\|\veps(\cdot,t)\|_{L^\infty(\Om)} \le \|v_0\|_{L^\infty(\Om)}
	\qquad \mbox{for all } t>0
  \ee
  as well as
  \be{uv}
	\int_0^\infty \io \frac{\ueps\veps}{1+\eps\ueps} \le \io v_0
	\qquad \mbox{for all } \eps\in (0,1).
  \ee
  Moreover, given any $T>0$ one can find $C(T)>0$ such that
  \be{lnv}
	\io \ln \frac{1}{\veps(\cdot,t)} \le C(T)
	\qquad \mbox{for all $t\in (0,T)$ and } \eps\in (0,1),
  \ee
  and that
  \be{grad_lnv}
	\int_0^T \io \frac{|\na\veps|^2}{\veps^2} \le C(T)
	\qquad \mbox{for all } \eps\in (0,1).
  \ee
\end{lem}
\proof
  This can be verified by an essentially verbatim copy of the arguments from 
  \cite[Lemma 2.2, Lemma 4.1]{win_2d}, supplemented by the observation that (\ref{uv}) holds 
  due to the fact that 
  \bas
	\int_0^t \io \frac{\ueps\veps}{1+\eps\ueps}
	= \io v_0 - \io \veps(\cdot,t) \le \io v_0
	\qquad \mbox{for all $t>0$ and } \eps\in (0,1)
  \eas
  according to the second equation in (\ref{0eps}).
\qed
Throughout the sequel, we shall henceforth fix a smoothly bounded $\Om\subset\R^n$ as well as initial data $(u_0,v_0)$ 
fulfilling (\ref{init}), and whenever a function $\phi$ satisfying (\ref{preg}) is given, we shall let 
$((\ueps,\veps))_{\eps\in (0,1)}$ denote the family of solutions to (\ref{0eps})-(\ref{pe}) accordingly obtained
in Lemma \ref{lem_loc}.\abs
A first common regularity feature of these solutions will result from a duality-based argument based on the following
observation which in its essence goes back to \cite{taowin_M3AS} already.
Here and below, for $\vp\in L^1(\Om)$ we abbreviate its average according to $\ov{\vp}:=\frac{1}{|\Om|} \io \vp$.
\begin{lem}\label{lem3}
  Let $D(A):=\{\vp\in W^{2,2}(\Om) \ | \ \io \vp=0 \mbox{ and } \frac{\pa\vp}{\pa\nu}=0 \mbox{ on } \pO \}$ and
  $A\vp:=-\Del\vp$ for $\vp\in D(A)$. Then whenever (\ref{preg}) holds, we have
  \be{3.1}
	\frac{1}{2} \frac{d}{dt} \io \big|A^{-\frac{1}{2}} (\ueps-\ouz)\big|^2
	+ \io \ueps^2 \pe(\veps)
	= \ouz^2 |\Om|\eps 	
	+ \ouz \io \ueps \phi(\veps)
	\qquad \mbox{for all $t>0$ and } \eps\in (0,1).
  \ee
\end{lem}
\proof
  Since $\io (\ueps-\ouz)=0$ for all $t>0$ and $\eps\in (0,1)$ by (\ref{mass}), and since clearly also
  $\io \big(\ueps\pe(\veps)-\ov{\ueps\pe(\veps)}\big)=0$ for all $t>0$ and $\eps\in (0,1)$, using that
  according to (\ref{0eps}) we have
  $\pa_t (\ueps-\ouz)=\Del\big(\ueps\pe(\veps)-\ov{\ueps\pe(\veps)}\big)$ in $\Om\times (0,\infty)$ for all $\eps\in (0,1)$, we see 
  that
  \bas
	\pa_t A^{-1} (\ueps-\ouz) 
	= - \ueps \pe(\veps) + \ov{\ueps\pe(\veps)}
	\quad \mbox{in } \Om\times (0,\infty)
	\qquad \mbox{for all } \eps\in (0,1).
  \eas
  We only need to multiply this by $\ueps-\ouz$ and use the self-adjointness of $A^{-\frac{1}{2}}$ to infer that, indeed,
  \bas
	\frac{1}{2} \frac{d}{dt} \io \big|A^{-\frac{1}{2}} (\ueps-\ouz)\big|^2
	&=& \io \big( -\ueps\pe(\veps) + \ov{\ueps\pe(\veps)}\big) (\ueps-\ouz) \\
	&=& - \io \ueps^2 \pe(\veps)
	+ \ouz \io \ueps\pe(\veps)
	+ \ov{\ueps\pe(\veps)} \io (\ueps-\ouz) \\
	&=& - \io \ueps^2 \pe(\veps)
	+ \ouz^2 |\Om| \eps
	+ \ouz \io \ueps \phi(\veps)
	\qquad \mbox{for all $t>0$ and } \eps\in (0,1)
  \eas
  due to (\ref{mass}). 
\qed
In a general setting compatible with the assumptions both of Theorem \ref{theo15} and Theorem \ref{theo155},
this can be seen to imply a first a priori estimate beyond those from Lemma \ref{lem_loc}. 
We announce already here, however, that in the context of nonlinearities $\phi$ which grow at most linearly
near the origin, Lemma \ref{lem44} will provide a significant refinement which will form the origin for
the time-independent bounds claimed in Theorem \ref{theo155}.
\begin{lem}\label{lem4}
  Suppose that (\ref{preg}) and (\ref{pa}) hold with some $\al>0$. 
  Then for all $T>0$ there exists $C(T)>0$ such that
  \be{4.1}
	\int_0^T \io \ueps^2 \pe(\veps) 
	+ \int_0^T \io \ueps^2 \veps^\al 
	\le C(T)
	\qquad \mbox{for all } \eps\in (0,1).
  \ee
\end{lem}
\proof
  Writing $c_1:=\|v_0\|_{L^\infty(\Om)}$ and $c_2:=\|\phi\|_{L^\infty((0,c_1))}$, from (\ref{vinfty}) we infer that
  $\phi(\veps)\le c_2$ in $\Om\times (0,\infty)$ for all $\eps\in (0,1)$, with finiteness of $c_2$ being guaranteed by (\ref{preg}).
  On the right-hand side of (\ref{3.1}), again using (\ref{mass}) we can accordingly estimate
  \bas
	\ouz^2 |\Om|\eps 	
	+ \ouz \io \ueps \phi(\veps)
	&\le& \ouz^2 |\Om|\eps 	
	+ c_2 \ouz \io \ueps
	= (\eps+c_2) \ouz^2 |\Om| \\
	&\le& (1+c_2) \ouz^2 |\Om|
	\qquad \mbox{for all $t>0$ and } \eps\in (0,1),
  \eas
  whence upon an integration in time we see that
  \be{4.2}
	\frac{1}{2} \io \big|A^{-\frac{1}{2}}(\ueps(\cdot,T)-\ouz)\big|^2
	+ \int_0^T \io \ueps^2 \pe(\veps)
	\le \frac{1}{2} \io \big|A^{-\frac{1}{2}}(u_0-\ouz)\big|^2
	+ (1+c_2) \ouz^2 |\Om| T
  \ee  
  for all $T>0$ and $\eps\in (0,1)$.
  Since (\ref{pa}) asserts positivity of $c_3:=\inf_{\xi\in (0,c_1)} \frac{\phi(\xi)}{\xi^\al}$, and since (\ref{vinfty}) ensures 
  that
  \bas
	\pe(\veps) \ge \frac{1}{2} \pe(\veps) + \frac{c_3}{2} \veps^\al
	\quad \mbox{in } \Om\times (0,\infty)
	\qquad \mbox{for all } \eps\in (0,1),
  \eas
  from (\ref{4.2}) we already obtain (\ref{4.1}).
\qed
Playing a key role in our reasoning, the standard logarithmic entropy can be described with respect to a very basic
evolution feature as follows.
\begin{lem}\label{lem21}
  If (\ref{preg}) holds, then
  \be{21.1}
	\frac{d}{dt} \io \ueps\ln\ueps
	+ \io \pe(\veps) \frac{|\na\ueps|^2}{\ueps}
	= -\io \pe'(\veps) \na\ueps\cdot\na\veps
	\qquad \mbox{for all $t>0$ and } \eps\in (0,1).
  \ee
\end{lem}
\proof
  This can directly be seen upon multiplying the first equation in (\ref{0eps}) by $\ln\ueps$ and integrating by parts using
  (\ref{mass}).
\qed
Thus led to provide appropriate control over the expression on the right of (\ref{21.1}) and especially the 
taxis gradients $\na\veps$, possibly with singular weights originating from potentially unbounded factors $\pe'(\veps)$,
we first perform a very standard testing procedure to obtain 
a general statement on a basic regularity property, which at this stage is yet conditional
by relying on a square integrablity feature of $(\ueps\veps)_{\eps\in (0,1)}$.
\begin{lem}\label{lem54}
  Assume (\ref{preg}). Then there exists $C>0$ such that
  \be{54.1}
	\io |\na\veps(\cdot,t)|^2 \le C
	+ \int_0^t \io \ueps^2\veps^2
	\qquad \mbox{for all $t>0$ and } \eps\in (0,1),
  \ee
  and that 
  \be{54.23}
	\int_0^T \io |\Del\veps|^2
	+ \int_0^T \io \frac{|\na\veps|^4}{\veps^2}
	+ \int_0^T \io v_{\eps t}^2 
	\le C + C \int_0^T \io \ueps^2\veps^2
	\qquad \mbox{for all $T>0$ and } \eps\in (0,1).
  \ee
\end{lem}
\proof
  We test the second equation in (\ref{0eps}) against $-\Del\veps$ and $v_{\eps t}$ in a standard manner to see that due to
  Young's inequality,
  \bea{54.4}
	\frac{d}{dt} \io |\na\veps|^2
	+ \frac{1}{2} \io |\Del\veps|^2
	+ \frac{1}{2} \io v_{\eps t}^2
	&=& - \frac{1}{2} \io |\Del\veps|^2
	- \frac{1}{2} \io v_{\eps t}^2
	+ \io \frac{\ueps\veps}{1+\eps\ueps} \Del\veps
	- \io \frac{\ueps\veps}{1+\eps\ueps} v_{\eps t} \nn\\
	&\le& \io \Big(\frac{\ueps\veps}{1+\eps\ueps}\Big)^2 \nn\\
	&\le& \io \ueps^2\veps^2
	\qquad \mbox{for all $t>0$ and } \eps\in (0,1).
  \eea
  Since
  \bas
	\io \frac{|\na\veps|^4}{\veps^2}
	&=& - \io |\na\veps|^2 \na\veps\cdot\na\frac{1}{\veps}
	= \io \frac{1}{\veps} \Big\{ \na\veps\cdot\na |\na\veps|^2 + |\na\veps|^2 \Del\veps\Big\} \\
	&\le& (2+\sqrt{n}) \io \frac{1}{\veps} |\na\veps|^2 |D^2\veps| \\
	&\le& \frac{1}{2} \io \frac{|\na\veps|^4}{\veps^2} + \frac{(2+\sqrt{n})^2}{2} \io |D^2\veps|^2
	\qquad \mbox{for all $t>0$ and } \eps\in (0,1),
  \eas
  and since thus
  \bas
	\io \frac{|\na\veps|^4}{\veps^2} \le (2+\sqrt{n})^2 \io |\Del\veps|^2
	\qquad \mbox{for all $t>0$ and } \eps\in (0,1)
  \eas
  due to the fact that according to the identity $\na\veps\cdot \na\Del\veps=\frac{1}{2}\Del |\na\veps|^2 - |D^2\veps|^2$ we have
  \bas
	\io |\Del\veps|^2
	= - \io \na\veps\cdot \na \Del\veps
	= \io |D^2\veps|^2 - \frac{1}{2} \int_{\pO} \frac{\pa |\na\veps|^2}{\pa\nu}
	\ge \io |D^2\veps|^2
	\quad \mbox{for all $t>0$ and } \eps\in (0,1)
  \eas
  thanks to the convexity of $\Om$ (\cite{lions_ARMA}), from (\ref{54.4}) we obtain both (\ref{54.1}) and
  (\ref{54.23}).
\qed
As a crucial preparation for our analysis of (\ref{21.1}) in the context of both Theorem \ref{theo15} and
the particular subcase $\al=1$ of Theorem \ref{theo155}, we note that
if integral bounds even for $\ueps^2\veps$ can be drawn on, then the taxis gradient can be controlled
even when weighted in a more singular manner than in (\ref{54.23}).
This can be confirmed in the course of another fairly well-established variational reasoning:
\begin{lem}\label{lem56}
  Assume (\ref{preg}). Then there exists $C>0$ such that
  \be{56.1}
	\int_0^T \io \frac{|\na\veps|^4}{\veps^3} 
	\le C + C \int_0^T \io \ueps^2 \veps
	\qquad \mbox{for all $T>0$ and } \eps\in (0,1).
  \ee
\end{lem}
\proof
  We integrate by parts using the second equation in (\ref{0eps}) to find that again since 
  $\na\veps\cdot\na\Del\veps=\frac{1}{2}\Del |\na\veps|^2 - |D^2\veps|^2$ for all $\eps\in (0,1)$,
  \bea{56.2}
	\frac{1}{2} \frac{d}{dt} \io \frac{|\na\veps|^2}{\veps}
	&=& \io \frac{\na\veps}{\veps} \cdot \Big\{ \na\Del\veps - \na\frac{\ueps\veps}{1+\eps\ueps}\Big\}
	- \frac{1}{2} \io \frac{|\na\veps|^2}{\veps^2} \cdot \Big\{ \Del\veps - \frac{\ueps\veps}{1+\eps\ueps}\Big\} \nn\\
	&=& \frac{1}{2} \io \frac{1}{\veps} \Del |\na\veps|^2
	- \io \frac{|D^2\veps|^2}{\veps}
	- \frac{1}{2} \io \frac{|\na\veps|^2}{\veps^2} \Del\veps \nn\\
	& & + \io \frac{\ueps\veps}{1+\eps\ueps} \na\cdot \Big( \frac{\na\veps}{\veps}\Big)
	+ \frac{1}{2} \io \frac{\ueps}{1+\eps\ueps} \frac{|\na\veps|^2}{\veps} \nn\\
	&=& - \io \frac{|D^2\veps|^2}{\veps}
	+ \io \frac{1}{\veps^2} \na\veps\cdot\na |\na\veps|^2
	- \io \frac{|\na\veps|^4}{\veps^3}
	+ \frac{1}{2} \int_{\pO} \frac{1}{\veps} \frac{\pa |\na\veps|^2}{\pa\nu} \nn\\
	& & + \io \frac{\ueps}{1+\eps\ueps} \Del\veps
	- \frac{1}{2} \io \frac{\ueps}{1+\eps\ueps} \frac{|\na\veps|^2}{\veps} \nn\\
	&\le& - \io \frac{|D^2\veps|^2}{\veps}
	+ \io \frac{1}{\veps^2} \na\veps\cdot\na |\na\veps|^2
	- \io \frac{|\na\veps|^4}{\veps^3} \nn\\
	& & + \io \frac{\ueps}{1+\eps\ueps} \Del\veps
	\qquad \mbox{for all $t>0$ and } \eps\in (0,1),
  \eea
  because $\frac{\pa |\na\veps|^2}{\pa\nu} \le 0$ on $\pO\times (0,\infty)$ for all $\eps\in (0,1)$ by convexity of
  $\Om$ (\cite{lions_ARMA}).
  Here we may use the well-known facts (\cite[p.331]{win_CPDE2012} and \cite[Lemma 3.4]{win_DCDSB2022}) that
  \bas
	- \io \frac{|D^2\veps|^2}{\veps}
	+ \io \frac{1}{\veps^2} \na\veps\cdot\na |\na\veps|^2
	- \io \frac{|\na\veps|^4}{\veps^3} 
	= - \io \veps |D^2 \ln \veps|^2
  \eas
  for all $t>0$ and $\eps\in (0,1)$, 
  and that with $c_1:=\frac{1}{2(3+\sqrt{n})}$ we have 
  \bas
	c_1 \io \frac{|\na\veps|^4}{\veps^3} + c_1 \io \frac{|D^2\veps|^2}{\veps} \le \io \veps |D^2\ln\veps|^2
	\qquad \mbox{for all $t>0$ and } \eps\in (0,1),
  \eas
  and employ Young's inequality to estimate
  \bas
	\io \frac{\ueps}{1+\eps\ueps} \Del\veps
	&\le& \sqrt{n} \io \ueps |D^2\veps| \\
	&\le& c_1 \io \frac{|D^2\veps|^2}{\veps}
	+ \frac{n}{4c_1} \io \ueps^2 \veps
	\qquad \mbox{for all $t>0$ and } \eps\in (0,1).
  \eas
  An integration of (\ref{56.2}) therefore shows that
  \bas
	\frac{1}{2} \io \frac{|\na\veps(\cdot,T)|^2}{\veps(\cdot,T)}
	+ c_1 \int_0^T \io \frac{|\na\veps|^4}{\veps^3}
	\le \frac{1}{2} \io \frac{|\na v_0|^2}{v_0}
	+ \frac{n}{4c_1} \int_0^T \io \ueps^2 \veps
	\qquad \mbox{for all $T>0$ and } \eps\in (0,1),
  \eas
  and hence establishes (\ref{56.1}) according to (\ref{init}).
\qed
The following analysis of the products $\ueps\veps^2$ will provide a handy
path toward the derivation of a favorable compactness feature which, thanks to the positivity feature of the $\veps$
encrypted in (\ref{lnv}), will form the source for pointwise a.e.~convergence of $(\ueps)_{\eps\in (0,1)}$ along
a suitable subsequence (cf.~Lemma \ref{lem14}).
\begin{lem}\label{lem33}
  Assume (\ref{preg}), (\ref{pa}) and (\ref{p2}) with some $\al\in (0,2]$, and let $k\in\N$ be such that $k>\frac{n+2}{2}$. 
  Then there exists $C>0$ such that
  \be{33.1}
	\io \big|\na (\ueps\veps^2)\big| 
	\le C \cdot \bigg\{ 
	\io \pe(\veps) \frac{|\na\ueps|^2}{\ueps}
	+ \io \frac{|\na\veps|^4}{\veps^2}
	+ \io \ueps^2\veps^2
	+1 \bigg\}
  \ee
  and
  \be{33.2}
	\big\| \pa_t (\ueps\veps^2)\big\|_{(W^{k,2}(\Om))^\star} 
	\le C \cdot \bigg\{ 
	\io \pe(\veps) \frac{|\na\ueps|^2}{\ueps}
	+ \io |\Del\veps|^2
	+ \io \frac{|\na\veps|^4}{\veps^2}
	+ \io \ueps^2\veps^2
	+1 \bigg\}
  \ee
  for all $t>0$ and $\eps\in (0,1)$.
\end{lem}
\proof
  Using (\ref{preg}), (\ref{pa}), (\ref{p2}) and (\ref{vinfty}), we fix positive constants $c_1, c_2, c_3$ and $c_4$ such that
  \be{33.3}
	\veps\le c_1,
	\quad
	c_2 \veps^\al \le \pe(\veps) \le c_3
	\quad \mbox{and} \quad
	|\pe'(\veps)|\le c_4 \veps^{\al-1}
	\quad \mbox{in $\Om\times (0,\infty)$ \quad for all } \eps\in (0,1),
  \ee
  whence due to Young's inequality and (\ref{mass}),
  \bea{33.4}
	\io \big|\na (\ueps\veps^2)\big|
	&\le& \io \veps^2 |\na\ueps|
	+ 2 \io \ueps\veps |\na\veps| \nn\\
	&\le& \frac{1}{2} \io \pe(\veps) \frac{|\na\ueps|^2}{\ueps}
	+ \frac{1}{2} \io \ueps \frac{\veps^4}{\pe(\veps)}
	+ \io \ueps^2 \veps^2
	+ \frac{1}{2} \io \frac{|\na\veps|^4}{\veps^2}
	+ \frac{1}{2} \io \veps^2 \nn\\
	&\le& \frac{1}{2} \io \pe(\veps) \frac{|\na\ueps|^2}{\ueps}
	+ \frac{c_1^{4-\al}}{2c_2} \ouz |\Om|
	+ \io \ueps^2 \veps^2
	+ \frac{1}{2} \io \frac{|\na\veps|^4}{\veps^2}
	+ \frac{c_1^2 |\Om|}{2}
  \eea
  for all $t>0$ an $\eps\in (0,1)$.
  Likewise, for fixed $\vp\in C^\infty(\bom)$ fulfilling $\|\vp\|_{L^\infty(\Om)} + \|\na\vp\|_{L^\infty(\Om)} \le 1$,
  recalling (\ref{0eps}) we see that
  \bea{33.5}
	\bigg| \io \pa_t (\ueps\veps^2) \vp \bigg|
	&=& \bigg| - \io \na (\ueps\pe(\veps)) \cdot \na (\veps^2\vp)
	+ 2\io \ueps \veps \cdot \Big\{ \Del \veps - \frac{\ueps\veps}{1+\eps\ueps}\Big\}\cdot \vp \bigg| \nn\\
	&=& \bigg| 
	-2 \io \veps \pe(\veps) (\na\ueps\cdot\na\veps) \vp
	- 2\io \ueps\veps \pe'(\veps) |\na\veps|^2 \vp \nn\\
	& & - \io \veps^2 \pe(\veps) \na\ueps\cdot\na\vp
	- \io \ueps \veps^2 \pe'(\veps) \na\veps\cdot\na\vp \nn\\
	& & + 2\io \ueps\veps \Del\veps \vp
	- 2\io \frac{\ueps^2\veps^2}{1+\eps\ueps} \vp \bigg| \nn\\
	&\le& 2\io \veps \pe(\veps) |\na\ueps| \cdot |\na\veps|
	+ 2c_4 \io \ueps\veps^\al |\na\veps|^2 \nn\\
	& & + \io \veps^2 \pe(\veps) |\na\ueps| 
	+ c_4 \io \ueps \veps^{\al+1} |\na\veps| \nn\\
	& & +2 \io \ueps\veps |\Del\veps| 
	+ 2\io \ueps^2\veps^2 \nn\\
	&\le& \io \pe(\veps) \frac{|\na\ueps|^2}{\ueps}
	+ \io \ueps\veps^2 \pe(\veps) |\na\veps|^2 \nn\\
	& & + 2c_4 \io \ueps\veps^\al |\na\veps|^2 \nn\\
	& & + \frac{1}{2} \io \pe(\veps) \frac{|\na\ueps|^2}{\ueps}
	+ \frac{1}{2} \io \ueps\veps^4 \pe(\veps) \nn\\
	& & + \frac{c_4}{2} \io \veps^{2\al} |\na\veps|^2 \nn\\
	& & + \io |\Del\veps|^2
	+ \Big(\frac{c_4}{2}+3\Big) \io \ueps^2\veps^2
	\qquad \mbox{for all $t>0$ and } \eps\in (0,1).
  \eea
  Here, again by Young's inequality, (\ref{33.3}) and (\ref{mass}),
  \bas
	& & \hs{-20mm}
	\io \ueps\veps^2 \pe(\veps) |\na\veps|^2 + 2c_4 \io \ueps\veps^\al |\na\veps|^2 
	+ \frac{c_4}{2} \io \veps^{2\al} |\na\veps|^2 \nn\\
	&\le& \frac{1}{2} \io \ueps^2\veps^2
	+ \frac{1}{2} \io \veps^2 \pe^2(\veps) |\na\veps|^4 \\
	& & + c_4 \io \ueps^2\veps^2
	+ c_4 \io \veps^{2\al-2}|\na\veps|^4 \nn\\
	& & + \frac{c_4}{4} \io \frac{|\na\veps|^4}{\veps^2}
	+ \frac{c_4}{4} \io \veps^{4\al+2} \nn\\
	&\le& \Big(c_4+\frac{1}{2}\Big) \io \ueps^2\veps^2
	+ \Big( \frac{c_1^4 c_3^2}{2} + c_1^{2\al} c_4 + \frac{c_4}{4}\Big) \io \frac{|\na\veps|^4}{\veps^2}
	+ \frac{c_1^{4\al+2} c_4 |\Om|}{4}
  \eas
  and
  \bas
	\frac{1}{2} \io \ueps \veps^4 \pe(\veps)
	\le \frac{c_1^4 c_3}{2} \io u_0
  \eas
  for all $t>0$ and $\eps\in (0,1)$.
  Therefore, (\ref{33.4}) implies (\ref{33.1}), whereas (\ref{33.2}) results from (\ref{33.5}), because $W^{k,2}(\Om) \hra
  W^{1,\infty}(\Om)$ due to our assumption that $k>\frac{n+2}{2}$.
\qed
\mysection{The case $\al\in (0,1)$. Proof of Theorem \ref{theo15}}\label{sect3}
In the particular setting specified in Theorem \ref{theo15}, a very rough estimation of the expressions
on the right-hand sides of (\ref{54.1}), (\ref{54.23}) and (\ref{56.1}) on the basis of Lemma \ref{lem4} immediately
yields time-dependent bounds for the second solution component in the following sense.
\begin{cor}\label{cor58}
  If (\ref{preg}) and (\ref{pa}) hold with some $\al\in (0,1)$, then for any $T>0$ there exists $C(T)>0$ such that
  \bas
	\io |\na\veps(\cdot,t)|^2 \le C(T)
	\qquad \mbox{for all $t\in (0,T)$ and } \eps\in (0,1),
  \eas
  and that
  \bas
	\int_0^T \io |\Del\veps|^2 + \int_0^T \io v_{\eps t}^2 
	+ \int_0^T \io \frac{|\na\veps|^4}{\veps^2} 
	+ \int_0^T \io \frac{|\na\veps|^4}{\veps^3}
	\le C(T)
	\qquad \mbox{for all } \eps\in (0,1).
  \eas
\end{cor}
\proof
  Since for all $T>0$ and $\eps\in (0,1)$ we have
  \bas
	\int_0^T \io \ueps^2\veps^2 \le \|v_0\|_{L^\infty(\Om)}^{2-\al} \int_0^T \io \ueps^2\veps^\al
	\quad \mbox{and} \quad
	\int_0^T \io \ueps^2\veps \le \|v_0\|_{L^\infty(\Om)}^{1-\al} \int_0^T \io \ueps^2\veps^\al
  \eas
  due to (\ref{vinfty}) and the inequality $\al<1$, this is a direct consequence of Lemma \ref{lem54},
  Lemma \ref{lem56} and Lemma \ref{lem4}.
\qed
The core of our analysis concerning Theorem \ref{theo15} can now be found launched in the following evolution feature
of $0\le t\mapsto \io \ueps(\cdot,t)\pe(\veps(\cdot,t))$.
\begin{lem}\label{lem10}
  Assume (\ref{preg}). Then
  \bea{10.1}
	\frac{d}{dt} \io \ueps\pe(\veps)
	&=& - \io \Big\{ \pe(\veps)\pe'(\veps)+\pe'(\veps)\Big\} \na\ueps\cdot\na\veps
	- \io \ueps \cdot \Big\{ \pe'^2(\veps) + \pe''(\veps)\Big\} |\na\veps|^2 \nn\\
	& & - \io \frac{\ueps^2\veps}{1+\eps\ueps} \pe'(\veps)
	\qquad \mbox{for all $t>0$ and } \eps\in (0,1).
  \eea
\end{lem}
\proof
  This follows from straightforward computation using (\ref{0eps}): Indeed, for all $t>0$ and $\eps\in (0,1)$ we have
  \bea{10.2}
	\frac{d}{dt} \io \ueps\pe(\veps)
	&=& \io \Del (\ueps\pe(\veps)) \cdot \pe(\veps)
	+ \io \ueps \pe'(\veps) \cdot \Big\{ \Del\veps - \frac{\ueps\veps}{1+\eps\ueps}\Big\} \nn\\
	&=& \io \ueps\pe(\veps) \Del\pe(\veps)
	+ \io \ueps\pe'(\veps) \Del\veps
	- \io \frac{\ueps^2\veps}{1+\eps\ueps} \pe'(\veps),
  \eea 
  where one further integration by parts shows that
  \bea{10.3}
	& & \hs{-25mm}
	\io \ueps\pe(\veps) \Del\pe(\veps)
	+ \io \ueps\pe'(\veps) \Del\veps \nn\\
	&=& \io \ueps\pe(\veps) \pe''(\veps) |\na\veps|^2
	+ \io \ueps \cdot \Big\{ \pe(\veps)\pe'(\veps) + \pe'(\veps)\Big\} \Del\veps \nn\\
	&=& \io \ueps\pe(\veps) \pe''(\veps) |\na\veps|^2
	- \io \Big\{ \pe(\veps)\pe'(\veps) + \pe'(\veps)\Big\} \na\ueps\cdot\na\veps \nn\\
	& & - \io \ueps\cdot \big\{ \pe\pe' + \pe' \big\}'(\veps) |\na\veps|^2
	\qquad \mbox{for all $t>0$ and } \eps\in (0,1).
  \eea
  Since
  \bas
	\pe\pe''- \big\{ \pe\pe' + \pe'\big\}'
	= - \pe'^2 - \pe''
	\quad \mbox{on $(0,\infty)$ \qquad for all } \eps\in (0,1),
  \eas
  from (\ref{10.2}) and (\ref{10.3}) we obtain (\ref{10.1}).
\qed
In fact, it will turn out that when appropriately combined with Lemma \ref{lem21}, the property (\ref{10.1})
will imply a quasi-energy feature of certain among the functionals in (\ref{energy1}).
Our selection of the number $a>0$ appearing therein will be accomplished in the course of the following elementary
but crucial analysis concerned with the behavior of $(\pe)_{\eps\in (0,1)}$ within finite intervals of the form
$[0,\xis]$, where $\xis$ will finally be chosen to coincide with $\|v_0\|_{L^\infty(\Om)}$.
\begin{lem}\label{lem11}
  Suppose that (\ref{preg}) and (\ref{p3}) hold with some $\al\in (0,1)$ and $\xi_0>0$.
  Then there exist $a>0$ and $\epss\in (0,1)$ with the property that whenever $\xis>0$, one can find $C(\xis)>0$ such that
  \be{11.1}
	\frac{\big\{ a\pe(\xi)\pe'(\xi) + (a-1) \pe'(\xi)\big\}^2}{2\pe(\xi)}
	+ a\pe'^2(\xi) + a\pe''(\xi) \le \frac{C(\xis)}{\xi}
	\qquad \mbox{for all $\xi\in [0,\xis]$ and } \eps\in (0,\epss).
  \ee
\end{lem}
\proof
  We let $a:=\frac{1}{\al}$ and then observe that since $\al<1$,
  \bas
	\lim_{\wt{\phi}\to 0} \Big\{ a^2 \wt{\phi}^2 +2a(a-1) \wt{\phi} + (a-1)^2 + 2a\wt{\phi} - \frac{2(1-\al)a}{\al} \Big\}
	&=& (a-1)^2 - \frac{2(1-\al)a}{\al}
	= a^2 - \frac{2}{\al} a + 1 \\
	&=& 1-\frac{1}{\al^2} 
	<0,
  \eas
  whence due to the identity $\phi(0)=0$ required in (\ref{preg}) we can fix $\xi_1\in (0,\xi_0)$ such that
  \bas
	\max_{\xi\in [0,\xi_1]} \Big\{ a^2 \phi^2(\xi) + 2a(a-1) \phi(\xi) + (a-1)^2 + 2a\phi(\xi) - \frac{2(1-\al)a}{\al}\Big\}
	<0.
  \eas
  Since clearly $\pe\to\phi$ in $L^\infty([0,\xi_1])$ as $\eps\searrow 0$, this entails the existence of $\epss\in (0,1)$ such that
  \be{11.2}
	a^2 \pe^2(\xi) + 2a(a-1) \pe(\xi) + (a-1)^2 + 2a\pe(\xi) - \frac{2(1-\al)a}{\al} \le 0
	\qquad \mbox{for all $\xi\in [0,\xi_1]$ and } \eps\in (0,\epss).
  \ee
  To derive (\ref{11.1}) from this, we note that thanks to (\ref{p3}),
  \bas
	0 \ge \big(\phi^\frac{1}{\al}\big)''(\xi)
	= \frac{1}{\al} \phi^{\frac{1}{\al}-1}(\xi)\phi''(\xi)
	+ \frac{1}{\al} \Big(\frac{1}{\al}-1\Big) \phi^{\frac{1}{\al}-2}(\xi)\phi'^2(\xi)
	\qquad \mbox{for all } \xi\in (0,\xi_0),
  \eas
  and that thus
  \bas
	\phi(\xi)\phi''(\xi) \le - \frac{1-\al}{\al} \phi'^2(\xi)
	\qquad \mbox{for all } \xi\in (0,\xi_0),
  \eas
  in particular implying that $\phi''\le 0$ on $(0,\xi_0)$ and hence also
  \bas
	\pe(\xi)\pe''(\xi) 
	&=& \phi(\xi)\phi''(\xi) + \eps\phi''(\xi)
	\le \phi(\xi)\phi''(\xi) \\
	&\le& - \frac{1-\al}{\al} \phi'^2(\xi) 
	= - \frac{1-\al}{\al} \pe'^2(\xi)
	\qquad \mbox{for all $\xi\in (0,\xi_0)$ and } \eps\in (0,1).
  \eas
  As $\xi_1\le\xi_0$, (\ref{11.2}) therefore guarantees that
  \bea{11.3}
	& & \hs{-20mm}
	\big\{ a\pe(\xi)\pe'(\xi)+(a-1) \pe'(\xi)\big\}^2
	+ 2a\pe(\xi)\pe'^2(\xi)
	+ 2a\pe(\xi)\pe''(\xi) \nn\\
	&\le& 	\big\{ a\pe(\xi)\pe'(\xi)+(a-1) \pe'(\xi)\big\}^2
	+ 2a\pe(\xi)\pe'^2(\xi) - \frac{2(1-\al)a}{\al} \pe'^2(\xi) \nn\\
	&=& \pe'^2(\xi) \cdot \Big\{ a^2 \pe^2(\xi) + 2a(a-1) \pe(\xi) + (a-1)^2 + 2a\pe(\xi) - \frac{2(1-\al)a}{\al}\Big\} 
		\nn\\[2mm]
	&\le& 0
	\qquad \mbox{for all $\xi\in [0,\xi_1]$ and } \eps\in (0,\epss).
  \eea
  Now given $\xis>0$, in the case when $\xis\le \xi_1$ we immediately infer (\ref{11.1}) from (\ref{11.3}).
  Otherwise, we use the regularity and positivity properties of $\phi$ asserted by (\ref{preg}) to see that, again since
  $\pe=\phi+\eps$ and thus $\pe'\equiv \phi'$ as well as $\pe''\equiv \phi''$ for all $\eps\in (0,1)$, we can find
  $c_1(\xis)>0$ such that
  \bas
	\frac{\big\{ a\pe(\xi)\pe'(\xi) + (a-1) \pe'(\xi)\big\}^2}{2\pe(\xi)}
	+ a\pe'^2(\xi) + a\pe''(\xi) \le c_1(\xis)
	\qquad \mbox{for all $\xi\in [\xi_1,\xis]$ and } \eps\in (0,1).
  \eas
  Once more in view of (\ref{11.3}), we thus infer that (\ref{11.1}) also holds in this case if we let $C(\xis):=\xis c_1(\xis)$.
\qed
Our main step toward Theorem \ref{theo15} can now be achieved by concatenating Lemma \ref{lem21} and Lemma \ref{lem10}
through Lemma \ref{lem11}.
\begin{lem}\label{lem12}
  Assume (\ref{preg}), (\ref{pa}), (\ref{p2}) and (\ref{p3}) with some $\al\in (0,1)$ and $\xi_0>0$.
  Then there exists $\epss\in (0,1)$ such that for each $T>0$ it is possible to fix $C(T)>0$ in such a way that
  \be{12.1}
	\io \ueps(\cdot,t) \ln \ueps(\cdot,t) \le C(T)
	\qquad \mbox{for all $t\in (0,T)$ and } \eps\in (0,\epss)
  \ee
  and
  \be{12.01}
	\int_0^T \io \pe(\veps) \frac{|\na\ueps|^2}{\ueps} \le C(T)
	\qquad \mbox{for all } \eps\in (0,\epss).
  \ee
\end{lem}
\proof
  We let $a>0$ and $\epss\in (0,1)$ be as provided by Lemma \ref{lem11}, and applying said lemma to $\xis:=\|v_0\|_{L^\infty(\Om)}$,
  thanks to (\ref{vinfty}) we can pick $c_1>0$ such that
  \be{12.2}
	\frac{\big\{ a\pe(\veps)\pe'(\veps)+(a-1)\pe'^2(\veps)\big\}^2}{2\pe(\veps)}
	+ a\pe'^2(\veps) + a\pe''(\veps) \le \frac{c_1}{\veps}
	\quad \mbox{in } \Om\times (0,\infty)
	\qquad \mbox{for all } \eps\in (0,\epss).
  \ee
  Keeping this value of $a$ fixed, we combine 
  Lemma \ref{lem21} with Lemma \ref{lem10} to see that
  \bea{12.3}
	& & \hs{-30mm}
	\frac{d}{dt} \bigg\{ \io \ueps\ln\ueps - a \io \ueps\pe(\veps)\bigg\} 
	+ \io \pe(\veps) \frac{|\na\ueps|^2}{\ueps} \nn\\
	&=& - \io \pe'(\veps) \na\ueps\cdot\na\veps \nn\\
	& & + a \io \Big\{ \pe(\veps)\pe'(\veps) + \pe'(\veps)\Big\} \na\ueps\cdot\na\veps \nn\\
	& & + a \io \ueps\cdot \Big\{ \pe'^2(\veps) + \pe''(\veps)\Big\} |\na\veps|^2 \nn\\
	& & + a \io \frac{\ueps^2 \veps}{1+\eps\ueps} \pe'(\veps)
	\qquad \mbox{for all $t>0$ and } \eps\in (0,1).
  \eea
  Here by Young's inequality and (\ref{12.2}),
  \bea{12.4}
	& & \hs{-14mm}
	- \io \pe'(\veps) \na\ueps\cdot\na\veps 
	+ a \io \Big\{ \pe(\veps)\pe'(\veps) + \pe'(\veps)\Big\} \na\ueps\cdot\na\veps 
	+ a \io \ueps\cdot \Big\{ \pe'^2(\veps) + \pe''(\veps)\Big\} |\na\veps|^2 \nn\\
	&=& \io \Big\{ a\pe(\veps)\pe'(\veps) + (a-1) \pe'(\veps)\Big\} \na\ueps\cdot\na\veps
	+ a \io \ueps\cdot \Big\{ \pe'^2(\veps) + \pe''(\veps)\Big\} |\na\veps|^2 \nn\\
	&\le& \frac{1}{2} \io \pe(\veps) \frac{|\na\ueps|^2}{\ueps}
	+ \io \ueps \cdot \bigg\{ \frac{\big\{ a\pe(\veps)\pe'(\veps)+(a-1)\pe'(\veps)\big\}^2}{2\pe(\veps)}
	+ a\pe'^2(\veps) + a\pe''(\veps) \bigg\} |\na\veps|^2 \nn\\
	&\le& \frac{1}{2} \io \pe(\veps) \frac{|\na\ueps|^2}{\ueps}
	+ c_1 \io \frac{\ueps}{\veps} |\na\veps|^2
	\qquad \mbox{for all $t>0$ and } \eps\in (0,\epss),
  \eea
  and again employing Young's inequality and using (\ref{vinfty}) we find that
  \bea{12.5}
	c_1 \io \frac{\ueps}{\veps} |\na\veps|^2
	&\le& \frac{c_1}{2} \io \frac{|\na\veps|^4}{\veps^3}
	+ \frac{c_1}{2} \io \ueps^2\veps \nn\\
	&\le& \frac{c_1}{2} \io \frac{|\na\veps|^4}{\veps^3}
	+ \frac{c_1 \|v_0\|_{L^\infty(\Om)}^{1-\al}}{2} \io \ueps^2\veps^\al 
	\qquad \mbox{for all $t>0$ and } \eps\in (0,1).
  \eea
  As, by (\ref{p2}) and (\ref{vinfty}), with some $c_2>0$ we have
  $|\pe'(\veps)| \le c_2 \veps^{\al-1}$ in $\Om\times (0,\infty)$ for all $\eps\in (0,1)$, we can furthermore estimate
  \bas
	a \io \frac{\ueps^2\veps}{1+\eps\ueps} \pe'(\veps)
	\le c_2 a \io \ueps^2 \veps^\al
	\qquad \mbox{for all $t>0$ and } \eps\in (0,1),
  \eas
  from (\ref{12.3}), (\ref{12.4}) and (\ref{12.5}) we infer that
  \bea{12.6}
	& & \hs{-20mm}
	\frac{d}{dt} \bigg\{ \io \ueps\ln\ueps - a \io \ueps\pe(\veps)\bigg\} 
	+ \frac{1}{2} \io \pe(\veps) \frac{|\na\ueps|^2}{\ueps} \nn\\
	&\le& \frac{c_1}{2} \io \frac{|\na\veps|^4}{\veps^3}
	+ \Big\{ \frac{c_1 \|v_0\|_{L^\infty(\Om)}^{1-\al}}{2} + c_2 a\Big\} \cdot \io \ueps^2 \veps^\al
	\qquad \mbox{for all $t>0$ and } \eps\in (0,\epss).
  \eea
  Since (\ref{vinfty}) together with (\ref{preg}) clearly entails the existence of $c_3>0$ such that
  \bas
	a\io \ueps\pe(\veps) \le c_3
	\qquad \mbox{for all $t>0$ and } \eps\in (0,1),
  \eas
  upon integrating (\ref{12.6}) we readily see that (\ref{12.1}) and (\ref{12.01}) are consequences of 
  Corollary \ref{cor58} and Lemma \ref{lem4}.
\qed
Based on the weighted estimate in (\ref{12.01}), we can additionally make sure that also the fluxes acting 
in the first equation from (\ref{0eps}) enjoy bounds in reflexive Lebesgue spaces determined by (\ref{15.4}).
\begin{lem}\label{lem13}
  Let (\ref{preg}), (\ref{pa}), (\ref{p2}) and (\ref{p3}) be satisfied with some $\al\in (0,1)$ and $\xi_0>0$.
  Then one can find $\epss\in (0,1)$ such that for each $T>0$ there exists $C(T)>0$ such that
  \be{13.1}
	\int_0^T \io \big| \na (\ueps\pe(\veps))\big|^{p(\al)} \le C(T)
	\qquad \mbox{for all } \eps\in (0,\epss),
  \ee
  where $p(\al)$ is as defined in (\ref{15.4}).
\end{lem}
\proof 
  According to (\ref{preg}), (\ref{p2}) and (\ref{vinfty}), there exist positive constants $c_1, c_2$ and $c_3$ such that
  \be{13.22}
	\veps\le c_1, \quad
	\pe(\veps)\le c_2
	\quad \mbox{and} \quad
	|\pe'(\veps)|\le c_3 \veps^{\al-1}
	\quad \mbox{in } \Om\times (0,\infty)
	\quad \mbox{for all } \eps\in (0,1),
  \ee
  whence fixing $T>0$ and $\eps\in (0,1)$ henceforth and writing $p=p(\al)$ we see that
  \bea{13.3}
	\int_0^T \io \big|\na (\ueps\pe(\veps))\big|^p
	&=& \int_0^T \io \big| \pe(\veps)\na\ueps + \ueps\pe'(\veps)\na\veps\big|^p \nn\\
	&\le& 2^{p-1} \int_0^T \io \pe^p(\veps) |\na\ueps|^p
	+ 2^{p-1} c_3^p \int_0^T \io \ueps^p \veps^{p(\al-1)} |\na\veps|^p.
  \eea
  Here since $p\le \frac{4}{3}$, we may use Young's inequality along with (\ref{13.22}) to estimate
  \bea{13.4}
	\int_0^T \io \pe^p(\veps) |\na\ueps|^p
	&\le& \int_0^T \io \pe^\frac{4}{3}(\veps) |\na\ueps|^\frac{4}{3}
	+ |\Om| T \nn\\
	&=& \int_0^T \io \Big( \pe(\veps) \frac{|\na\ueps|^2}{\ueps} \Big)^\frac{2}{3} 
		\cdot \ueps^\frac{2}{3} \pe^\frac{2}{3}(\veps)
	+ |\Om| T \nn\\
	&\le& \int_0^T \io \pe(\veps) \frac{|\na\ueps|^2}{\ueps}
	+ \int_0^T \io \ueps^2 \pe^2(\veps)
	+ |\Om| T \nn\\
	&\le& \int_0^T \io \pe(\veps) \frac{|\na\ueps|^2}{\ueps}
	+ c_2 \int_0^T \io \ueps^2 \pe(\veps)
	+ |\Om| T,
  \eea
  and in the case when $\al\ge\frac{1}{2}$, we can proceed similarly to find that
  \bea{13.5}
	\int_0^T \io \ueps^p \veps^{p(\al-1)} |\na\veps|^p
	&=& \int_0^T \io \Big( \frac{|\na\veps|^4}{\veps^3} \Big)^\frac{1}{3} \cdot \ueps^\frac{4}{3} \veps^\frac{4\al-1}{3} \nn\\
	&\le& \int_0^T \io \frac{|\na\veps|^4}{\veps^3}
	+ \int_0^T \io \ueps^2 \veps^\frac{4\al-1}{2} \nn\\
	&\le& \int_0^T \io \frac{|\na\veps|^4}{\veps^3}
	+ c_1^\frac{2\al-1}{2} \int_0^T \io \ueps^2 \veps^\al.
  \eea
  If $\al<\frac{1}{2}$, then again by Young's inequality,
  \bea{13.6}
	\int_0^T \io \ueps^p \veps^{p(\al-1)} |\na\veps|^p
	&=& \int_0^T \io \big( \ueps^2 \veps^\al \big)^\frac{p}{2} \veps^\frac{p(\al-2)}{2} |\na\veps|^p \nn\\
	&\le& \int_0^T \io \ueps^2 \veps^\al
	+ \int_0^T \io \veps^\frac{p(\al-2)}{2-p} |\na\veps|^\frac{2p}{2-p},
  \eea
  where according to (\ref{15.4}) we then have $\frac{p}{2-p}=\frac{1}{1-\al}$ and thus, once more by Young's inequality,
  \bea{13.7}
	\int_0^T \io \veps^\frac{p(\al-2)}{2-p} |\na\veps|^\frac{2p}{2-p}
	&=& \int_0^T \io \veps^\frac{\al-2}{1-\al} |\na\veps|^\frac{2}{1-\al} \nn\\
	&=& \int_0^T \io \Big( \frac{|\na\veps|^2}{\veps^2} \Big)^\frac{1-2\al}{1-\al} \cdot
		\veps^{-\frac{3\al}{1-\al}} |\na\veps|^\frac{4\al}{1-\al} \nn\\
	&\le& \int_0^T \io \frac{|\na\veps|^2}{\veps^2}
	+ \int_0^T \io \frac{|\na\veps|^4}{\veps^3}.
  \eea
  In view of (\ref{13.3})-(\ref{13.7}), the claim thus becomes a consequence of Lemma \ref{lem12}, Lemma \ref{lem4},
  Corollary \ref{cor58} and (\ref{grad_lnv}).
\qed
A straightforward subsequence extraction now yields a global solution in the flavor of Theorem \ref{theo15}.
\begin{lem}\label{lem14}
  Assume (\ref{preg}), (\ref{pa}), (\ref{p2}) and (\ref{p3}) with some $\al\in (0,1)$ and $\xi_0>0$, and let (\ref{init}) hold.
  Then there exist $(\eps_j)_{j\in\N} \subset (0,1)$ as well as functions $u$ and $v$ fulfilling (\ref{reg})-(\ref{15.3})
  such that $\eps_j\searrow 0$ as $j\to\infty$, that $u\ge 0$ and $v>0$ a.e.~in $\Om\times (0,\infty)$,
  and that 
  \begin{eqnarray}
	& & \ueps \to u
	\qquad \mbox{in $L^1_{loc}(\bom\times [0,\infty))$ and a.e.~in $\Om\times (0,\infty)$,}
	\label{14.2} \\
	& & \veps \to v
	\qquad \mbox{in $L^p_{loc}(\bom\times [0,\infty))$ for all $p\in [1,\infty)$ and a.e.~in } \Omega\times (0,\infty)
	\label{14.3}
	\qquad \mbox{and} \\
	& & \na\veps\to \na v
	\qquad \mbox{in } L^2_{loc}(\bom\times [0,\infty))
	\label{14.4}
  \end{eqnarray}
  as $\eps=\eps_j\searrow 0$.
  Moreover, (\ref{massu}) holds, and $(u,v)$ is a global weak solution of (\ref{0}) in the sense of Definition \ref{dw}.
\end{lem}
\proof
  For fixed $T>0$, from Lemma \ref{lem12} we particularly obtain that with $\epss>0$ as introduced there,
  \be{14.5}
	(\ueps)_{\eps\in (0,\epss)} 
	\mbox{ is uniformly integrable over $\Om\times (0,T)$,} 
  \ee
  while using Corollary \ref{cor58} we find that
  \be{14.6}
	(\veps)_{\eps\in (0,1)}
	\mbox{ is bounded in } L^\infty((0,T);W^{1,2}(\Om)) \mbox{ and in } L^2((0,T);W^{2,2}(\Om)),
  \ee
  and that
  \be{14.7}
	(v_{\eps t})_{\eps\in (0,1)}
	\mbox{ is bounded in } L^2(\Om\times (0,T)),
  \ee
  which in line with an Aubin-Lions lemma (\cite{temam}) particularly means that
  \be{14.8}
	(\veps)_{\eps\in (0,1)}
	\mbox{ is relatively compact with respect to the strong topology in } L^2((0,T);W^{1,2}(\Om)).
  \ee
  Apart from that, Lemma \ref{lem33} together with Lemma \ref{lem12}, Corollary \ref{cor58}, Lemma \ref{lem4},
  (\ref{mass}) and (\ref{vinfty}) ensures that
  \bas
	(\ueps\veps^2)_{\eps\in (0,1)}
	\mbox{ is bounded in } L^1((0,T);W^{1,1}(\Om))
  \eas
  and that whenever $k\in\N$ is such that $k>\frac{n+2}{2}$,
  \bas
	\big(\pa_t(\ueps\veps^2)\big)_{\eps\in (0,1)}
	\mbox{ is bounded in } L^1\big( (0,T);(W^{k,2}(\Om)^\star\big),
  \eas
  and that thus, once more by an Aubin-Lions lemma,
  \bas
	(\ueps\veps^2)_{\eps\in (0,1)}
	\mbox{ is relatively compact in } L^1(\Om\times (0,T)).
  \eas
  In view of a standard extraction argument using (\ref{vinfty}), we accordingly infer the existence of
  $(\eps_j)_{j\in\N}\subset (0,\epss)$ such that $\eps_j\searrow 0$ as $j\to\infty$, and that (\ref{14.3}) and (\ref{14.4}) 
  as well as
  \be{14.111}
	\ueps\veps^2 \to z
	\quad \mbox{a.e.~in } \Om\times (0,\infty)
	\qquad \mbox{as } \eps=\eps_j\searrow 0
  \ee
  with some nonnegative functions $v\in L^\infty(\Om\times (0,\infty)) \cap L^\infty_{loc}([0,\infty);W^{1,2}(\Om))
  \cap L^2_{loc}([0,\infty);W^{2,2}(\Om))$ and $z\in L^1_{loc}(\bom\times [0,\infty))$, about which due to (\ref{lnv})
  and Fatou's lemma we even know that $v>0$ a.e.~in $\Om\times (0,\infty)$.
  Therefore, letting $u:=\frac{z}{v^2}$ defines an a.e.~in $\Om\times (0,\infty)$ finite measurable function 
  $u$ for which thanks to (\ref{14.111}) and (\ref{14.3}) we have $\ueps\to u$ a.e.~in $\Om\times (0,\infty)$
  as $\eps=\eps_j\searrow 0$, and which thus satisfies (\ref{14.2}), and hence due to (\ref{mass}) also (\ref{massu}), 
  as a consequence of (\ref{14.5}) and the Vitali
  convergence theorem.
  The inequality in (\ref{15.2}) hence results from (\ref{14.2}) and Lemma \ref{lem12} upon an application of Fatou's lemma,
  and to deduce (\ref{15.3}) we note that since (\ref{14.5}), (\ref{preg}) and (\ref{vinfty}) moreover entail that 
  also $(\ueps \pe(\veps))_{\eps\in (0,\epss)}$ is uniformly integrable over $\Om\times (0,T)$, by means of (\ref{14.3}) 
  and again the Vitali convergence theorem we infer that furthermore 
  \be{14.10}
	\ueps\pe(\veps) \to u\phi(v) \quad \mbox{in } L^1_{loc}(\bom\times [0,\infty))
  \ee
  as $\eps=\eps_j\searrow 0$.
  This enables us to identify corresponding weak limits obtained upon employing Lemma \ref{lem13}, according to which, namely,
  we know that with $p(\al)>1$ taken from (\ref{15.4}), $(\na (\ueps\pe(\veps)))_{\eps\in (0,1)}$ is bounded, and hence 
  relatively compact with respect to the weak topology, in $L^p(\Om\times (0,T))$ for all $T>0$; 
  therefore, (\ref{14.10}) implies that 
  \be{14.11}
	\na\big(\ueps\pe(\veps)\big) \wto \na\big(u\phi(v)\big) \quad \mbox{in } L^{p(\al)}_{loc}(\bom\times [0,\infty))
  \ee
  as $\eps=\eps_j\searrow 0$, and that thus $(u,v)$ especially satisfies (\ref{15.3}).\abs
  Finally, for arbitrary $\vp\in C_0^\infty(\bom\times [0,\infty))$, the identities in (\ref{wu}) and (\ref{wv}) can be derived
  in a straightforward manner from (\ref{14.2})-(\ref{14.4}) and (\ref{14.10}), 
  relying on the fact that according to (\ref{14.2}) and (\ref{14.3}), and once more thanks to the Vitali theorem,
  also $\frac{\ueps\veps}{1+\eps\ueps} \to uv$ in
  $L^1_{loc}(\bom\times [0,\infty))$ as $\eps=\eps_j\searrow 0$.
\qed
Our main result on global solvability in the weakly degenerate case has, in fact, thereby been established:\abs
\proofc of Theorem \ref{theo15}. \quad
  The statement has fully been covered by Lemma \ref{lem14} already.
\qed
\mysection{The case $\al\in [1,2]$. Proof of Theorem \ref{theo155}}
In the context of the assumptions from Theorem \ref{theo155}, essentially due to local Lipschitz continuity of 
$\phi$ thereby implied we can refine our analysis already at a rather early stage as follows.
\begin{lem}\label{lem44}
  Assume (\ref{preg}), (\ref{pa}) and (\ref{p2}) with some $\al\ge 1$.
  Then there exists $C>0$ such that
  \be{44.1}
	\int_0^T \io \ueps^2 \pe(\veps) 
	+ \int_0^T \io \ueps \veps^\al 
	\le C \Ieps(T)
	\qquad \mbox{for all $T>0$ and } \eps\in (0,1),
  \ee
  where we have set
  \be{I}
	\Ieps(T):=1+\eps T + \int_0^T \io \frac{\eps\ueps^2 \veps}{1+\eps\ueps}
	\qquad \mbox{for $T>0$ and } \eps\in (0,1).
  \ee
\end{lem}
\proof
  Let $c_1:=\|v_0\|_{L^\infty(\Om)}$. Then noting that $c_2:=\|\phi'\|_{L^\infty([0,c_1])}$ is finite by (\ref{p2}) and
  our assumption that $\alpha\ge 1$, using (\ref{vinfty}) we see that $\pe(\veps)\le c_2 \veps+\eps$ in $\Om\times (0,\infty)$
  for all $\eps\in (0,1)$. Therefore, an integration of (\ref{3.1}) shows that
  \bea{44.2}
	\frac{1}{2} \io \big|A^{-\frac{1}{2}}(\ueps(\cdot,T)-\ouz)\big|^2
	+ \int_0^T \io \ueps^2 \pe(\veps)
	&\le& \frac{1}{2} \io \big|A^{-\frac{1}{2}}(u_0-\ouz)\big|^2
	+ 2 \ouz^2 |\Om|\eps T 
	+ c_2 \ouz \int_0^T \io \ueps\veps \nn\\[1mm]
	& & \hs{20mm}
	\qquad \mbox{for all $T>0$ and } \eps\in (0,1).
  \eea
  Here, again relying on (\ref{vinfty}) we may estimate $\pe(\veps) \ge c_3\veps^\al$ in $\Om\times (0,\infty)$ for all
  $\eps\in (0,1)$, with $c_3:=\inf_{\xi\in [0,c_1]} \frac{\phi(\xi)}{\xi^\al}$ being positive by (\ref{pa}), so that again
  \bas
	\int_0^T \io \ueps^2 \pe(\veps) 
	\ge \frac{1}{2} \int_0^T \io \ueps^2 \pe(\veps) + \frac{c_3}{2} \int_0^T \io \ueps^2 \veps^\al
	\qquad \mbox{for all $T>0$ and } \eps\in (0,1).
  \eas
  Since furthermore
  \bas
	\int_0^T \io \ueps\veps
	&=& \int_0^T \io \frac{\ueps\veps}{1+\eps\ueps} 
	+ \int_0^T \io \frac{\eps \ueps^2\veps}{1+\eps\ueps} \\
	&\le& \io v_0
	+ \int_0^T \io \frac{\eps \ueps^2\veps}{1+\eps\ueps}
	\qquad \mbox{for all $T>0$ and } \eps\in (0,1)
  \eas
  by (\ref{uv}), from (\ref{44.2}) we thus obtain (\ref{44.1}) with some appropriately large $C>0$.
\qed
Now of crucial relevance in an appropriate handling of the stronger degeneracies from Theorem \ref{theo155}
will be the observation that not only for $\al\in (0,1)$ as in the previous part, but also for some $\al\ge 1$, 
the accordingly weaker regularity information then gained from bounds for $\int_0^T \io \ueps^2\veps^\al$ 
can be turned into expedient knowledge on the respective second solution components.
Our considerations in this direction will be rooted in the following outcome of elementary calculus.
\begin{lem}\label{lem60}
  Let $\al\in (1,2)$. Then for all $\vp\in C^2(\bom)$ fulfilling $\vp>0$ in $\bom$ and $\frac{\pa\vp}{\pa\nu}=0$ on $\pO$, we have
  \bea{60.1}
	& & \hs{-20mm}
	-2\io \vp^{\al-2} |D^2\vp|^2
	-(\al-2) \io \vp^{\al-3} \na\vp \cdot\na |\na \vp|^2
	+ (\al-2) \io \vp^{\al-3} |\na\vp|^2 \Del\vp \nn\\
	&=& - \frac{2}{(\al-1)^2} \io \vp^{-\al+2} |D^2 \vp^{\al-1}|^2
	- (\al-1)(2-\al) \io \vp^{\al-4} |\na\vp|^4
  \eea
  as well as
  \be{60.2}
	\io \vp^{\al-2} |D^2\vp|^2
	\le \frac{2}{(\al-1)^2} \io \vp^{-\al+2} |D^2 \vp^{\al-1}|^2
	+ 2(\al-2)^2 \io \vp^{\al-4} |\na\vp|^4
  \ee
  for any such $\vp$.
\end{lem}
\proof
  In view of a standard approximation procedure, we only need to consider the case when additionally
  $\vp\in C^3(\bom)$. 
  We may then integrate by parts to rewrite
  \bas
	& & \hs{-20mm}
	-(\al-2) \io \vp^{\al-3} \na\vp \cdot\na |\na \vp|^2
	+ (\al-2) \io \vp^{\al-3} |\na\vp|^2 \Del\vp \nn\\
	&=& - 4(\al-2) \io \vp^{\al-3} \na\vp\cdot (D^2\vp\cdot\na\vp)
	- (\al-2)(\al-3) \io \vp^{\al-4} |\na\vp|^4,
  \eas
  and to thus obtain (\ref{60.1}) by using the pointwise identity
  \be{60.3}
	|D^2\vp|^2 = \frac{1}{(\al-1)^2} \vp^{4-2\al} |D^2 \vp^{\al-1}|^2
	- 2(\al-2) \frac{1}{\vp} \na\vp \cdot (D^2\vp\cdot\na\vp)
	- (\al-2)^2 \frac{|\na\vp|^4}{\vp^2}.
  \ee
  Since, again by (\ref{60.3}), and by Young's inequality, any such $\vp$ satisfies
  \bas
	|D^2\vp|^2
	\le \frac{1}{(\al-1)^2} \vp^{4-2\al} |D^2\vp^{\al-1}|^2
	+ \bigg\{ \frac{1}{2} |D^2 \vp|^2
	+ 2(\al-2)^2 \frac{|\na \vp|^4}{\vp^2} \bigg\}
	- (\al-2)^2 \frac{|\na\vp|^4}{\vp^2}
  \eas
  and thus
  \bas
	|D^2\vp|^2
	\le \frac{2}{(\al-1)^2} \vp^{4-2\al} |D^2\vp^{\al-1}|^2
	+ 2(\al-2)^2 \frac{|\na \vp|^4}{\vp^2},
  \eas
  the inequality in (\ref{60.2}) follows after multiplication by $\vp^{\al-2}$ and integration.
\qed
When (\ref{preg}) and (\ref{pa}) hold with some $\al\in (1,2)$, 
controlling cross-diffusive gradients through (\ref{44.1})
thereby becomes possible by means of another independent testing procedure
applied to the second equation from (\ref{0eps}):
\begin{lem}\label{lem6}
  Assume (\ref{preg}) and (\ref{pa}) with some $\al\in (1,2)$.
  Then there exists $C>0$ such that
  \be{6.2}
	\int_0^T \io \veps^{\al-4} |\na\veps|^4 \le C \Ieps(T)
	\qquad \mbox{for all $T>0$ and } \eps\in (0,1),
  \ee
  where $(\Ieps)_{\eps\in (0,1)}$ is as defined through (\ref{I}).
\end{lem}
\proof
  We use the second equation in (\ref{0eps}) and integrate by parts to see that again since for all $\eps\in (0,1)$ we have
  $\na\veps\cdot\na\Del\veps=\frac{1}{2}\Del |\na\veps|^2 - |D^2\veps|^2$ in $\Om\times (0,\infty)$ and
  $\frac{\pa |\na\veps|^2}{\pa\nu} \le 0$ on $\pO\times (0,\infty)$ by convexity of $\Om$ (\cite{lions_ARMA}),
  \bea{6.3}
	\frac{d}{dt} \io \veps^{\al-2} |\na\veps|^2
	&=& 2 \io \veps^{\al-2} \cdot \Big\{ \na\veps\cdot\na\Del\veps - \na\veps\cdot \na \frac{\ueps\veps}{1+\eps\ueps} \Big\} 
		\nn\\
	& & + (\al-2) \io \veps^{\al-3} |\na\veps|^2 \cdot \Big\{ \Del\veps-\frac{\ueps\veps}{1+\eps\ueps}\Big\} \nn\\
	&=& \io \veps^{\al-2} \Del |\na\veps|^2
	- 2 \io \veps^{\al-2} |D^2\veps|^2
	+ (\al-2) \io \veps^{\al-3} |\na\veps|^2 \Del\veps \nn\\
	& & + 2\io \frac{\ueps\veps}{1+\eps\ueps} \na\cdot (\veps^{\al-2} \na\veps) 
	- (\al-2) \io \frac{\ueps\veps^{\al-2}}{1+\eps\ueps} |\na\veps|^2 \nn\\
	&=& -(\al-2) \io \veps^{\al-3} \na\veps\cdot\na |\na\veps|^2
	+ \int_{\pO} \veps^{\al-2} \frac{\pa |\na\veps|^2}{\pa\nu} \nn\\
	& & - 2 \io \veps^{\al-2} |D^2\veps|^2
	+ (\al-2) \io \veps^{\al-3} |\na\veps|^2 \Del\veps \nn\\
	& & + 2\io \frac{\ueps\veps^{\al-1}}{1+\eps\ueps} \Del\veps
	+ (\al-2) \io \frac{\ueps\veps^{\al-2}}{1+\eps\ueps} |\na\veps|^2 \nn\\
	&\le& -(\al-2) \io \veps^{\al-3} \na\veps\cdot\na |\na\veps|^2
	- 2 \io \veps^{\al-2} |D^2\veps|^2
	+ (\al-2) \io \veps^{\al-3} |\na\veps|^2 \Del\veps \nn\\
	& & + 2\io \frac{\ueps\veps^{\al-1}}{1+\eps\ueps} \Del\veps
	\qquad \mbox{for all $t>0$ and } \eps\in (0,1),
  \eea
  because $\al\le 2$. Here according to (\ref{60.1}),
  \bea{6.4}
	& & \hs{-30mm}
	-(\al-2) \io \veps^{\al-3} \na\veps\cdot\na |\na\veps|^2
	- 2 \io \veps^{\al-2} |D^2\veps|^2
	+ (\al-2) \io \veps^{\al-3} |\na\veps|^2 \Del\veps \nn\\
	&=& - \frac{2}{(\al-1)^2} \io \veps^{-\al+2} |D^2\veps^{\al-1}|^2
	- (\al-1)(2-\al) \io \veps^{\al-4} |\na\veps|^4
  \eea
  for all $t>0$ and $\eps\in (0,1)$,
  while thanks to the second statement in Lemma \ref{lem60} we can fix $c_1>0$ such that 
  \bea{6.5}
	\frac{2}{(\al-1)^2} \io \veps^{-\al+2} |D^2\veps^{\al-1}|^2
	+\frac{(\al-1)(2-\al)}{2} \io \veps^{\al-4} |\na\veps|^4
	\ge c_1 \io \veps^{\al-2} |\Del\veps|^2
  \eea
  for all $t>0$ and $\eps\in (0,1)$.
  We thereupon employ Young's inequality to see that
  \bas
	2\io \frac{\ueps\veps^{\al-1}}{1+\eps\ueps} \Del\veps
	&\le& c_1 \io \veps^{\al-2} |\Del\veps|^2
	+ \frac{1}{c_1} \io \frac{\ueps^2\veps^\al}{(1+\eps\ueps)^2} \\
	&\le& c_1 \io \veps^{\al-2} |\Del\veps|^2
	+ \frac{1}{c_1} \io \ueps^2 \veps^\al
	\qquad \mbox{for all $t>0$ and } \eps\in (0,1),
  \eas
  whence combining this with (\ref{6.3})-(\ref{6.5}) shows that
  \bas
	\frac{d}{dt} \io \veps^{\al-2} |\na\veps|^2
	+ \frac{(\al-1)(2-\al)}{2} \io \veps^{\al-4} |\na\veps|^4
	\le \frac{1}{c_1} \io \ueps^2 \veps^\al
	\qquad \mbox{for all $t>0$ and } \eps\in (0,1).
  \eas
  In view of Lemma \ref{lem4}, an integration in time completes the proof.
\qed
Combining this with previous knowledge yields the following summary of regularity features enjoyed by $(\veps)_{\eps\in (0,1)}$
throughout the entire range of parameters $\al\in [1,2]$.
\begin{cor}\label{cor57}
  If (\ref{preg}) and (\ref{pa}) hold with some $\al\in [1,2]$, then there exists $C>0$ such that with $(I_\eps)_{\eps\in (0,1)}$
  taken from (\ref{I}) we have
  \be{57.01}
	\io |\na\veps(\cdot,t)|^2 \le C\Ieps(t)
	\qquad \mbox{for all $t\in (0,T)$, any $T>0$ and each } \eps\in (0,1)
  \ee
  as well as
  \be{57.02}
	\int_0^T \io |\Del\veps|^2
	+ \int_0^T \io v_{\eps t}^2 \le C\Ieps(T)
	\qquad \mbox{for all $T>0$ and } \eps\in (0,1)
  \ee
  and
  \be{57.1}
	\int_0^T \io \veps^{\al-4} |\na\veps|^4 \le C\Ieps(T)
	\qquad \mbox{for all $T>0$ and } \eps\in (0,1).
  \ee
\end{cor}
\proof
  The estimates in (\ref{57.01}) and (\ref{57.02}) directly result from Lemma \ref{lem54} and Lemma \ref{lem4}.
  To confirm (\ref{57.1}), we note that in the case when $\al=1$, this immediately follows from Lemma \ref{lem56} 
  when combined with Lemma \ref{lem4},
  while if $\al=2$, then we may similarly conclude using Lemma \ref{lem54} together with Lemma \ref{lem4}.
  When $\al\in (1,2)$, (\ref{57.1}) has precisely been asserted by Lemma \ref{lem6}.
\qed
Instead of relying on (\ref{energy1}), we can now directly estimate the integral on the right-hand side of (\ref{21.1})
to achieve the following boundedness information which partially parallels that from Lemma \ref{lem12} in its outcome, but 
which significantly differs from the latter with regard to its derivation.
\begin{lem}\label{lem7}
  Suppose that there exists $\al\in [1,2]$ such that
  (\ref{preg}), (\ref{pa}) and (\ref{p2}) hold.
  Then there exists $C>0$ such that with $(\Ieps)_{\eps\in (0,1)}$ taken from (\ref{I}) we have
  \be{7.1}
	\io \ueps(\cdot,t) \ln \ueps(\cdot,t) \le C \Ieps(T)
	\qquad \mbox{for all $t\in (0,T)$, each $T>0$ and any } \eps\in (0,1)
  \ee
  as well as
  \be{7.01}
	\int_0^T \io \pe(\veps) \frac{|\na\ueps|^2}{\ueps} \le C \Ieps(T)
	\qquad \mbox{for all $T>0$ and } \eps\in (0,1).
  \ee
\end{lem}
\proof
  Based on (\ref{p2}), (\ref{pa}) and (\ref{vinfty}), we fix $c_1>0$ and $c_2>0$ such that $|\pe'(\veps)|\le c_1\veps^{\al-1}$ and
  $\pe(\veps)\ge c_2\veps^\al$ in $\Om\times (0,\infty)$ for all $\eps\in (0,1)$.
  Then relying on Young's inequality, from Lemma \ref{lem21} we infer that
  \bas
	\frac{d}{dt} \io \ueps\ln\ueps
	+ \io \pe(\veps) \frac{|\na\ueps|^2}{\ueps}
	&=& -\io \pe'(\veps) \na\ueps\cdot\na\veps \\
	&\le& \frac{1}{2} \io \pe(\veps) \frac{|\na\ueps|^2}{\ueps} 
	+ \frac{1}{2} \io \ueps \frac{\pe'^2(\veps)}{\pe(\veps)} |\na\veps|^2 \\
	&\le& \frac{1}{2} \io \pe(\veps) \frac{|\na\ueps|^2}{\ueps} 
	+ \frac{c_1^2}{2c_2} \io \ueps \veps^{\al-2} |\na\veps|^2 \\
	&\le& \frac{1}{2} \io \pe(\veps) \frac{|\na\ueps|^2}{\ueps} 
	+ \frac{c_1^2}{4c_2} \io \veps^{\al-4} |\na\veps|^4 \\
	& & + \frac{c_1^2}{4c_2} \io \ueps^2 \veps^\al
	\qquad \mbox{for all $t>0$ and } \eps\in (0,1).
  \eas
  The claim therefore results upon an integration using Lemma \ref{lem4} and Corollary \ref{cor57}.
\qed
Again, integrability features of fluxes can be gained by suitable interpolation:
\begin{lem}\label{lem133}
  Assume (\ref{preg}), (\ref{pa}) and (\ref{p2}) with some $\al\in [1,2]$.
  Then there exists $C>0$ such that
  \be{133.1}
	\int_0^T \io \big| \na (\ueps\pe(\veps))\big|^\frac{4}{3} \le C\Ieps(T)
	\qquad \mbox{for all $T>0$ and } \eps\in (0,1),
  \ee
  where again $(\Ieps)_{\eps\in (0,1)}$ is taken from (\ref{I}).
\end{lem}
\proof 
  We proceed similarly as in Lemma \ref{lem13} to see that thanks to Young's inequality, 
  (\ref{preg}), (\ref{p2}) and (\ref{vinfty}), we can find positive constants $c_1$ and $c_2$ such that for
  any $t>0$ and $\eps\in (0,1)$ we have
  \bea{133.2}
	\hs{-8mm}
	\io \big| \na (\ueps\pe(\veps))\big|^\frac{4}{3} 
	&\le& c_1 \io \pe^\frac{4}{3}(\veps) |\na\ueps|^\frac{4}{3}
	+ c_1 \io \ueps^\frac{4}{3} \veps^\frac{4(\al-1)}{3} |\na\veps|^\frac{4}{3} \nn\\
	&\le& c_1 \io \pe(\veps) \frac{|\na\ueps|^2}{\ueps}
	+ c_1 \io \ueps^2 \pe^2(\veps)
	+ c_1 \io \veps^{\al-4} |\na\veps|^4
	+ c_1 \io \ueps^2 \veps^\frac{3\al}{2} \nn\\
	&\le& c_1 \io \pe(\veps) \frac{|\na\ueps|^2}{\ueps}
	+ c_2 \io \ueps^2 \pe(\veps)
	+ c_1 \io \veps^{\al-4} |\na\veps|^4
	+ c_2 \io \ueps^2 \veps^\al.
  \eea
  In view of Lemma \ref{lem7}, Lemma \ref{lem4} and Corollary \ref{cor57} an integration of (\ref{133.2}) 
  yields (\ref{133.1}).
\qed
Now to prepare an exploitation of the estimates from 
Lemma \ref{lem7}, Corollary \ref{cor57} and Lemma \ref{lem133} on fixed time intervals, we record a rough preliminary
estimate for the expressions in (\ref{I}).
\begin{lem}\label{lem132}
  If (\ref{preg}), (\ref{pa}) and (\ref{p2}) hold with some $\al\ge 1$, then the numbers from (\ref{I}) satisfy
  \be{132.1}
	\Ieps(T)\le 1+T + \ouz \|v_0\|_{L^\infty(\Om)} |\Om| T
	\qquad \mbox{for all $T>0$ and } \eps\in (0,1).
  \ee
\end{lem}
\proof
  Trivially estimating $\frac{\eps\ueps}{1+\eps\ueps} \le 1$ for $\eps\in (0,1)$, from (\ref{vinfty}) and (\ref{mass}) 
  we obtain that
  \bas
	\io \frac{\eps\ueps^2 \veps}{1+\eps\ueps}
	\le \io \ueps\veps
	\le \ouz \|v_0\|_{L^\infty(\Om)} |\Om|
	\qquad \mbox{for all $t>0$ and } \eps\in (0,1),
  \eas
  whence (\ref{132.1}) directly results from (\ref{I}).
\qed
Having this at hand, we can adapt the reasoning from Lemma \ref{lem14}, supplemented by an additional argument 
deriving (\ref{155.1}) and (\ref{155.2}) by means of a refined treatment of $(\Ieps)_{\eps\in (0,1)}$, to construct
global solutions in the sense claimed in Theorem \ref{theo155}:
\begin{lem}\label{lem144}
  Suppose that there exists $\al\in [1,2]$ such that (\ref{preg}), (\ref{pa}) and (\ref{p2}) hold, and assume (\ref{init}).
  Then there exists $(\eps_j)_{j\in\N} \subset (0,1)$ 
  such that $\eps_j\searrow 0$ as $j\to\infty$ and that (\ref{14.2})-(\ref{14.4}) hold with some functions
  $u$ and $v$ which satisfy (\ref{reg}) with $u\ge 0$ and $v>0$ a.e.~in $\Om\times (0,\infty)$,
  for which (\ref{massu}), (\ref{155.1}) and (\ref{155.2}) hold, and which are such that $(u,v)$ forms a global weak solution
  of (\ref{0}) in the sense of Definition \ref{dw}.
\end{lem}
\proof
  In a first step estimating the elements of $(\Ieps)_{\eps\in (0,1)}$ via Lemma \ref{lem132}, we can proceed in much the same 
  manner as in Lemma \ref{lem14}, this time relying on 
  Lemma \ref{lem7}, Corollary \ref{cor57} and Lemma \ref{lem133},
  to extract $(\eps_j)_{j\in\N}$ such that $\eps_j\searrow 0$ as $j\to\infty$, and that 
  (\ref{14.2})-(\ref{14.4}) are valid with some global weak solution $(u,v)$ of (\ref{0}) which satisfy (\ref{reg}) as well as
  $u\ge 0$ and $v>0$ a.e.~in $\Om\times (0,\infty)$, and for which (\ref{massu}) holds.
  It thus remains to verify the additional regularity properties in (\ref{155.1}) and (\ref{155.2}), 
  for which purpose we may now rely
  on (\ref{14.2}) and (\ref{14.3}) to see that $\frac{\eps\ueps^2\veps}{1+\eps\ueps} \to 0$ a.e.~in $\Om\times (0,\infty)$
  as $\eps=\eps_j\searrow 0$, so that since
  \bas
	0 \le \frac{\eps\ueps^2\veps}{1+\eps\ueps}
	= \frac{\eps\ueps}{1+\eps\ueps} \veps \ueps
	\le \|v_0\|_{L^\infty(\Om)} \ueps
	\quad \mbox{in } \Om\times (0,\infty)
	\qquad \mbox{for all } \eps\in (0,1)
  \eas
  according to (\ref{vinfty}), the $L^1$ approximation feature in (\ref{14.2}) ensures that for each fixed $T>0$,
  \bas
	\int_0^T \io \frac{\eps\ueps^2\veps}{1+\eps\ueps} \to 0
	\qquad \mbox{as } \eps=\eps_j\searrow 0
  \eas
  due to the dominated convergence theorem. Hence, for any such $T$ we infer from (\ref{I}) that
  \bas
	\Ieps(T) \to 1
	\qquad \mbox{as } \eps=\eps_j\searrow 0,
  \eas
  so that revisiting 
  Lemma \ref{lem7}, Corollary \ref{cor57} and Lemma \ref{lem133} provides $c_1>0$, $c_2>0$ and $c_3>0$ with the property that
  whenever $T>0$, we can find $\eps_0(T)\in (0,1)$ such that
  \bas
	\io \ueps\ln\ueps \le c_1
	\quad \mbox{and} \quad
	\io |\na\veps|^2\le c_2
	\qquad \mbox{for all $t\in (0,T)$ and } \eps\in (\eps_j)_{j\in\N} \cap (0,\eps_0(T))
  \eas
  as well as
  \bas
	\int_0^T \io \veps^{\al-4} |\na\veps|^4 
	+ \int_0^T \io \big|\na (\ueps\pe(\veps))\big|^\frac{4}{3} \le c_3
	\qquad \mbox{for all } \eps\in (\eps_j)_{j\in\N} \cap (0,\eps_0(T)).
  \eas
  In view of (\ref{14.2})-(\ref{14.4}), Fatou's lemma and lower semicontinuity of $L^p$ norms with respect to weak convergence,
  taking $\eps=\eps_j\searrow 0$ we therefore readily obtain that
  \bas
	\io u\ln u \le c_1
	\quad \mbox{and} \quad
	\io |\na v|^2 \le c_2
	\qquad \mbox{for a.e.~} t>0,
  \eas
  and that
  \bas
	\int_0^T \io v^{\al-4} |\na v|^4
	+ \int_0^T \io \big| \na (u\phi(v))\big|^\frac{4}{3} \le c_3
	\qquad \mbox{for all } T>0,
  \eas
  meaning that indeed both (\ref{155.1}) and (\ref{155.2}) hold.
\qed
Also in the more strongly degenerate setting addressed in Theorem \ref{theo155}, we have thus found global solutions
which even enjoy the additional boundedness and decay features expressed in (\ref{155.1}) and (\ref{155.2}):\abs
\proofc of Theorem \ref{theo155}.\quad
  We only need to apply Lemma \ref{lem144}.
\qed

\bigskip

{\bf Acknowledgement.} \quad
The author acknowledges support of the {\em Deutsche Forschungsgemeinschaft} (Project No.~462888149).


\begin{thebibliography}{99}
%
\bibitem{ahn_yoon}
  \sc Ahn, J., Yoon, C.: 
  \it Global well-posedness and stability of constant equilibria in parabolic-elliptic chemotaxis systems without gradient sensing.
  \rm Nonlinearity {\bf 32}, 1327-1351 (2019)
\bibitem{amann}
  \sc Amann, H.:
  \it Dynamic theory of quasilinear parabolic systems III. Global existence.
  \rm Math. Z. {\bf 202}, 219-250 (1989)
\bibitem{maini_et_al_JTB2013}
  \sc Belmonte-Beitia, J., Woolley, T.E., Scott, J.G., Maini, P.K., Gaffney, E.A.:
  \it Modelling biological invasions: Individual to population scales at interfaces.
  \rm J.~Theor.~Biol. {\bf334}, 1-12 (2013)
\bibitem{burger}
  \sc Burger, M., Lauren\c{c}ot, Ph., Trescases, A.: \it Delayed blow-up for chemotaxis models with local sensing.
  \rm J. London Math. Soc. {\bf 103}, 1596-1617 (2021)
\bibitem{engwer_hunt_surulescu}
  \sc Engwer, C., Hunt, A., Surulescu, C.: \it Effective equations for anisotropic glioma spread with proliferation: 
  a multiscale approach and comparisons with previous settings.
  \rm Math.~Med.~Biol. {\bf 33}, 435-459 (2016)
\bibitem{fu}
  \sc Fu, X., Tang, L.H., Liu, C., Huang, J.D., Hwa, T., Lenz, P.:
  \it Stripe formation in bacterial systems with density-suppressed motility.
  \rm Phys. Rev. Lett. {\bf 108}, 198102 (2012)
\bibitem{fujie_jiang_JDE2020}
  \sc Fujie, K., Jiang, J.: \it Global existence for a kinetic model of pattern formation with density-suppressed motilities
  \rm J.~Differential Equations {\bf 269}, 5338-5378 (2020)
\bibitem{fujie_jiang_arxiv}
  \sc Fujie, K., Jiang, J.: 
  \it Comparison methods for a Keller-Segel-type model of pattern formations with density-suppressed motilities.
  \rm Calc. Var. Partial Differential Equations {\bf 60}, 92 (2021)
\bibitem{fujie_jiang_ACAP2021}
  \sc Fujie, K., Jiang, J.: 
  \it Boundedness of classical solutions to a degenerate Keller-Segel type model with signal-dependent motilities.
  \rm Acta Appl.~Math. {\bf 176}, 3 (2021)
\bibitem{fujie_senba}
  \sc Fujie, K., Senba, T.: \it Global existence and infinite time blow-up of classical solutions 
  to chemotaxis systems of local sensing in higher dimensions.
  \rm Nonlinear Anal., Theory Methods Appl. {\bf 222}, 112987 (2022)
\bibitem{jiang_arxiv}
  \sc Jiang, J.: \it Boundedness and exponential stabilization in a parabolic-elliptic Keller-Segel model with 
  signal-dependent motilities for local sensing chemotaxis.
  \rm Acta Math. Sci., Ser. B, Engl. Ed. {\bf 42}, 825-846 (2022)
\bibitem{jiang_laurencot}
  \sc Jiang, J., Lauren\c{c}ot, Ph.: \it Global existence and uniform boundedness in a chemotaxis model 
  with signal-dependent motility.
  \rm J. Differential Equations {\bf 299}, 513-541 (2021)
\bibitem{jin_kim_wang}
  \sc Jin, H.-Y., Kim, Y.-J., Wang, Z.-A.:
  \it Boundedness, stabilization, and pattern formation driven by density-suppressed motility.
  \rm SIAM J.~Appl.~Math. {\bf 78}, 1632-1657 (2018)
\bibitem{kawasaki_JTB1997}
  \sc Kawasaki, K., Mochizusi, A., Matsushita, M., Umeda, T., Shigesada, N.:
  \it Modeling Spatio-Temporal Patterns Generated by Bacillus subtilis.
  \rm J. Theor. Biol. {\bf 188}, 177-185 (1997)
\bibitem{KS1}
  \sc Keller, E.F., Segel, L.A.: \it Initiation of slime mold aggregation viewed as an instability. 
  \rm J. Theor. Biol. {\bf 26}, 399-415 (1970)
\bibitem{KS3}
  \sc Keller, E.F., Segel, L.A.: \it Model for chemotaxis.
  \rm J. Theoret. Biol. {\bf 30}, 225-234 (1971)
\bibitem{laurencot_CMS}
  \sc Lauren\c{c}ot, Ph.: 
  \it  Long term spatial homogeneity for a chemotaxis model with local sensing and consumption.
  \rm Commun. Math. Sci. {\bf 21}, 1743-1750 (2023)
\bibitem{plaza}
  \sc Leyva, J.F., M\'alaga, C., Plaza, R.G.: 
  \it The effects of nutrient chemotaxis on bacterial aggregation patterns with non-linear degenerate cross diffusion.
  \rm Physica A {\bf 392}, 5644-5662 (2013)
\bibitem{liwin2} 
  \sc Li, G., Winkler, M.: \it Relaxation in a Keller-Segel-consumption system involving signal-dependent motilities.
  \rm Commun. Math. Sci. {\bf 21}, 299-322 (2023)
\bibitem{li_zhao_ZAMP}
  \sc Li, D., Zhao, J.: \it Global boundedness and large time behavior of solutions to a chemotaxis-consumption system 
  with signal-dependent motility.
  \rm Z.~Angew.~Math.~Phys., {\bf 72}, 57 (2021)
\bibitem{lions_ARMA}
  \sc Lions, P.L.: \it R\'esolution de probl\`emes elliptiques quasilin\'eaires.
  \rm Arch.~Rat.~Mech.~Anal. {\bf 74}, 335-353 (1980)
\bibitem{liu}
  \sc Liu, C., et al.:
  \it Sequential establishment of stripe patterns in an expanding cell population.
  \rm Science {\bf 334}, 238 (2011)
\bibitem{wenbin_lv_ZAMP}
  \sc Lv, W., Wang, Q.: 
  \it Global existence for a class of chemotaxis systems with signal-dependent motility, indirect signal production 
  and generalized logistic source.
  \rm Z.~Angew.~Math.~Physik {\bf 71}, 53 (2020)
\bibitem{wenbin_lv_EECT}
  \sc Lv, W., Wang, Q.: 
  \it Global existence for a class of Keller-Segel model with signal-dependent motility and general logistic term.
  \rm Evol. Equ. Control Theory. {\bf 10}, 25-36 (2021)
\bibitem{wenbin_lv_PROCA}
  \sc Lv, W., Wang, Q.: 
  \it A $n$-dimensional chemotaxis system with signal-dependent motility and generalized logistic source: Global existence 
  and asymptotic stabilization.
  \rm Proc. Roy. Soc. Edinburgh Sect. A {\bf 151}, 821-841 (2021)
\bibitem{othmer_stevens}
  \sc Othmer, H.G., Stevens, A.: \it Aggregation, blowup and collapse: The ABC's of taxis in reinforced random walks. 
  \rm SIAM J. Appl. Math. {\bf 57}, 1044-1081 (1997)
\bibitem{taowin_M3AS}
  \sc Tao, Y., Winkler, M.: 
  \it Effects of signal-dependent motilities in a Keller-Segel-type reaction-diffusion system 
  \rm Math.~Mod.~Meth.~Appl.~Sci. {\bf 27}, 1645-1683 (2017)
\bibitem{temam}
  \sc Temam, R.: \it Navier-Stokes equations. Theory and numerical analysis.
  Studies in Mathematics and its Applications. Vol. 2. 
  \rm North-Holland, Amsterdam, 1977
\bibitem{wang_wang}
  \sc Wang, J., Wang, M.: 
  \it Boundedness in the higher-dimensional Keller-Segel model with signal-dependent motility and logistic growth.
  \rm J.~Math.~Phys. {\bf 60}, 011507 (2019)
\bibitem{win_CPDE2012}
  \sc Winkler, M.: \it Global large-data solutions in a chemotaxis-(Navier-)Stokes system modeling cellular swimming in fluid drops.
  \rm Commun.~Partial Differential Equations {\bf 37}, 319-351 (2012)
\bibitem{win_NON}
  \sc Winkler, M.: 
  \it Can simultaneous density-determined enhancement of diffusion and cross-diffusion foster boundedness in Keller-Segel type 
  systems involving signal-depen\-dent motilities?
  \rm Nonlinearity {\bf 33}, 6590-6623 (2020)
\bibitem{win_DCDSB2022}
  \sc Winkler, M.: \it Approaching logarithmic singularities in quasilinear chemotaxis-consumption systems with
  signal-dependent sensitivities.
  \rm Discrete Contin.~Dyn.~Syst.~Ser.~B {\bf 27}, 6565-6587 (2022)
\bibitem{win_NON2023}
  \sc Winkler, M.:
  \it  Stabilization despite pervasive strong cross-degeneracies in a nonlinear diffusion model 
  for migration-consumption interaction.
  \rm Nonlinearity {\bf 36}, 4438-4469 (2023)
\bibitem{win_2d}
  \sc Winkler, M.: \it Application of the Moser-Trudinger inequality in the constuction of global solutions 
  to a strongly degenerate migration model.
  \rm Bull. Math. Sci. {\bf 13}, 2250012 (2023)
\bibitem{win_sig_dep_mot_SMP}
  \sc Winkler, M.: \it A quantitative strong parabolic maximum principle 
  and application to a taxis-type migration-consumption model involving signal-dependent degenerate diffusion.
  \rm Ann.~Inst.~H.~Poincar\'e Anal.~Non Lin\'eaire, to appear
\bibitem{yifu_wang}
  \sc Xu, C., Wang, Y.: \it Asymptotic behavior of a quasilinear Keller-Segel system with signal-suppressed motility.
  \rm Calc. Var. Partial Differential Equations {\bf 60}, 183 (2021)
%
\end{thebibliography}
\end{document}